\documentclass[a4,11pt]{article}
\usepackage[centertags]{amsmath}
\usepackage{amsmath,amsfonts,amssymb,amsthm}

\usepackage[latin1]{inputenc}
\newcommand{\bb}{\begin{equation}}
\newcommand{\ee}{\end{equation}}
\newcommand{\ds}{\displaystyle }

\def\p{\partial}

\def\ln{\mbox{\rm ln}}

\def\x1{{\xi }_{xx}}
\def\x2{{\xi }_{yy}}
\def\x3{{\xi }_{xy}}

\def\e1{{\eta }_{xx}}
\def\e2{{\eta }_{yy}}
\def\e3{{\eta }_{xy}}

\def\Ga{\Gamma}

\def\de{\delta}

\def\l1{{\lambda}_1}

\def\nb{\nabla}
\def\kd{\partial}

\def\lb{\Delta_{g}}
\def\f{\frac}
\def\G{\Gamma}
\def\bb{\begin{equation}}
\def\ee{\end{equation}}
\def\ba{\begin{array}}
\def\ea{\end{array}}

\def\R{\mathbb{R}}
\def\nb{\nabla}
\def\L{{\cal L}}

\begin{document}
\pagenumbering{arabic}
\title{\huge \bf Special Conformal Groups of a
Riemannian Manifold and Lie Point Symmetries of the Nonlinear
Poisson Equation}
\author{\rm \large Yuri Bozhkov$^1$ and Igor Leite Freire$^{1,2}$\\ \\
\it $^1$Instituto de Matem\'atica,
Estat\'\i stica e \\ \it Computa\c c\~ao Cient\'\i fica - IMECC \\
\it Universidade Estadual de Campinas - UNICAMP \\ \it C.P.
$6065$, $13083$-$970$ - Campinas - SP, Brasil
\\ \rm E-mail: bozhkov@ime.unicamp.br \\ \ \ \ \ \
igor@ime.unicamp.br \\ \\
\it $^2$Centro de Matemática, Computação e Cognição\\ \it
Universidade Federal do ABC - UFABC\\ \it Rua Catequese, $242$,
Jardim,
$09090-400$\\\it Santo André - SP, Brasil\\
\rm E-mail: igor.freire@ufabc.edu.br}
\date{\ }
\maketitle \vspace{1cm}
\begin{abstract}
We obtain a complete group classification of the Lie point
symmetries of nonlinear Poisson equations on generic (pseudo)
Riemannian manifolds $M$. Using this result we study their Noether
symmetries and establish the respective conservation laws. It is
shown that the projection of the Lie point symmetries on $M$ are
special subgroups of the conformal group of $M$. In particular, if
the scalar curvature of $M$ vanishes, the projection on $M$ of the
Lie point symmetry group of the Poisson equation with critical
nonlinearity is the conformal group of the manifold. We illustrate
our results by applying them to the Thurston geometries.
\end{abstract}
\vskip 1cm
\begin{center}
{2000 AMS Mathematics Classification numbers:\vspace{0.2cm}\\
35J50, 35J20, 35J60

\

 Key words: Lie point symmetry, Noether symmetry, conservation
laws, conformal group}
\end{center}

\newpage

\section{Introduction}

\

The study of differential equations on manifolds is the corner
stone of the Geometric Analysis. For this purpose various methods
have been applied: fixed point theorems, continuity method,
maximum principles, a priori estimates, Schauder theory, etc.
However it seems that it is little known how the symmetries of the
considered equation (or system) and the geometry of the manifold
are related. (By a `symmetry' we understand a Lie point symmetry
\cite{ba,bk,i,ol,ole,ov}.)

Let $M^n$ be a (pseudo) Riemannian manifold of dimension $n\geq 3$
endowed with a (pseudo) Riemannian metric $g=(g_{ij})$ given in
local
 coordinates $\{x^{1},\cdots,x^{n}\}$. In this
paper we shall study the Lie point symmetries of the Poisson
equation on $M^n$:
 \bb\label{1.1} \lb u+f(u)=0, \ee where
$$\lb u=\f{1}{\sqrt{g}}\f{\p}{\p x^{i}}(\sqrt{g}g^{ij}\f{\p u}{\p x^{j}})=
g^{ij}\nabla_{i}\nabla_{j}u=\nabla^{j}\nabla_{j}u=\nabla_{i}\nabla^{i}u$$
is the Laplace-Beltrami operator, $f$ is a smooth function,
$(g^{ij})$ is the inverse matrix of $(g_{ij})$, ${\nabla }_i$ is
the covariant derivative corresponding to the Levi-Civita
connection and we have used the Einstein summation convention,
that is, summation from $1$ to $n$ over repeated indices is
understood.

The equation (\ref{1.1}) can be equivalently written as \bb
\label{a3} H\equiv g^{ij}u_{ij}-\G^{i}u_{i}+f(u)=0,\ee where
$\Gamma^{i}:=g^{pq}\G_{pq}^{i}$, $\G_{pq}^{i} $ being the
Christoffel symbols, $u_i=\frac{\kd u}{\kd x^i}$, $u_{ij}=
\frac{{\kd }^2 u}{\kd x^i \kd x^j}$.

We observe that equation (\ref{1.1}) includes elliptic and
hyperbolic equations, depending on whether $(M^{n},g)$ is a
Riemannian or a pseudo Riemannian manifold. Such equations appear
in various geometric and mathematical physics contexts which we
shall not going into here. We merely mention the Poisson equations
in $\R^{n}$, in particular those involving critical exponents,
taking in (\ref{1.1}) $f(u)=u^{\f{n+2}{n-2}}$ and the Euclidean
metric; the Klein-Gordon equation, taking  in (\ref{1.1}) the
metric $ds^{2}=-dt^{2}+dx^{2}+dy^{2}+dz^{2}$ and $f(u)=u$; the
semilinear wave equations in $\R^{1+n}$, with $f''(u)\neq 0$ and
the metric $ds^{2}=-dt^{2}+\delta_{ij}dx^{i}dx^{j}$ in
(\ref{1.1}); and the Klein-Gordon equation on the $\mathbb{S}^{2}$
sphere. These particular equations have been studied in \cite{i,
ol, bm, bh, fu1, sv, A, igor}. The interested reader may also
consult the book \cite{i} of Ibragimov, where various aspects of
symmetry analysis of differential equations on manifolds are
presented. We would also like to mention the paper by Ratto and
Rigoli \cite{rr} in which these authors obtain gradient bounds and
Liouville type theorems for the Poisson equation (\ref{1.1}) on
complete Riemannian manifolds.

The purpose of this paper is three-fold. First we shall obtain a
complete group classification for the semilinear Poisson equations
(\ref{1.1}) by applying the S. Lie symmetry theory. Then we shall
study the Noether symmetries of (\ref{1.1}). (Noether symmetry =
variational or divergence symmetry.) The latter will be used to
establish the corresponding conservation laws via the Noether's
theorem.

Since we suppose that the reader is familiar with the basic
notions and methods of contemporary group analysis
\cite{ba,bk,i,ol,ole,ov}, we shall not present preliminaries
concerning Lie point symmetries of differential equations. For a
geometric viewpoint of Lie point symmetries, see
\cite{olvjdg,man}. We would just like to recall that, following
Olver (\cite{ol}, p. 182), to perform a group classification on a
differential equation
 involving a generic function $f$ consists of finding the Lie point symmetries of the
 given equation with arbitrary $f$, and, then, to determine all possible particular
forms of $f$ for which the symmetry group
 can be enlarged. Usually there
 exists a geometrical or physical motivation for considering such specific cases.

 Our first result is the following group
 classification theorem.

 \

 {\bf Theorem 1.} \rm 1.) {\it The Lie point symmetry group of the
Poisson equation $(\ref{1.1})$ with an arbitrary $f(u)$ coincides
with the isometry group of $M^{n}$. In this case the infinitesimal
 generator is given by}
 \bb\label{z1} X=\xi =\xi^{i}(x)\f{\p}{\p x^{i}} \ee
 {\it where
 \[{\cal L}_{\xi}g_{ij}=0 \]
 and $ {\cal L}_{\xi}$ is the Lie derivative with respect to the
 vector filed $\xi $.

 For some special choices of the function $f(u)$ it can be extended
  in the cases listed below.}

 \rm 2.) {\it If $f(u)=0$ then the symmetries have the form
 \bb\label{z2}X=\xi^{i}(x)\f{\p}{\p x^{i}} +[(\frac{2-n}{4}\;\mu (x)
 +c)\; u +b(x)]\frac{\p}{\p u}, \ee
 where $c$ is an arbitrary constant,
 \bb\label{z3} \lb b =0,\ee
 \bb\label{z4} \lb \mu =0\ee
 and $\xi =\xi^{i}(x)\ds{\f{\p}{\p x^{i}}}$ is a conformal Killing
 vector field such that}
 \bb\label{z5} {\cal L}_{\xi}g_{ij}=\mu g_{ij}.\ee

 \rm 3.) {\it In the case $f(u)=k=const\neq 0$ the symmetries are generated
 by}
 \bb\label{h1} X=\xi^{i}(x)\f{\p}{\p x^{i}}
 +[(\frac{n-2}{n+2}(\f{1}{k} (\lb b )-c)\;u +
 b(x)]\frac{\p}{\p u},\ee
 {\it where $c$ is an arbitrary constant, \[ {\lb }^{\!\!\!2} b =0, \] and
 $\xi =\xi^{i}(x)\ds{\f{\p}{\p x^{i}}}$ is a conformal Killing
 vector field such that}
 \bb\label{h2} {\cal L}_{\xi}g_{ij}=\f{4}{(n+2)}(c-\f{1}{k}\lb b) g_{ij}.\ee

 \rm 4.) {\it If the function $f$ is a linear function: $f(u)=u$,
 then the symmetry generator is given by $(\ref{z2})$ with $\xi =\xi^{i}(x)\f{\p}{\p
 x^{i}}$ satisfying} $(\ref{z5})$,
 \bb\label{z6} \lb b +b=0, \ee
 \bb\label{z7} \f{2-n}{4}\lb \mu +\mu =0. \ee

 \rm 5.) {\it For exponential nonlinearity $f(u)=e^u$ we have the
 generator
 \bb\label{z8} X=\xi^{i}(x)\f{\p}{\p x^{i}} - \mu  \frac{\p}{\p
 u},\ee
 where $\mu $ is a constant and $\xi =\xi^{i}(x)\f{\p}{\p x^{i}}$ is a homothety such
 that}
 \bb\label{z9} {\cal L}_{\xi}g_{ij}=\mu g_{ij},\;\;\;\mu
 =const.\ee

 \rm 6.) {\it For power nonlinearity $f(u)=u^p$, $p\neq 0$, $p\neq 1$,
 the Lie point symmetry group is generated by
 \bb\label{z10} X=\xi^{i}(x)\f{\p}{\p x^{i}} +\frac{ \mu }{1-p}\;u\; \frac{\p}{\p
 u},\ee
 where $\mu $ is a constant and the vector field $\xi =\xi^{i}(x)\f{\p}{\p x^{i}}$ is a
 homothety such that}
 \bb\label{z11} {\cal L}_{\xi}g_{ij}=\mu g_{ij},\;\;\;\mu
 =const.\ee

 \rm 6.1) {\it If $p=\ds{\frac{n+2}{n-2}}$, $n\neq 6$, the infinitesimal
 generator of the Lie point symmetries has the form
 \bb\label{z12} X=\xi^{i}(x)\f{\p}{\p x^{i}} +\frac{ 2-n }{4}\mu \;u\; \frac{\p}{\p
 u},\ee
 where $\mu $ is a harmonic function on $M^n$:
 \bb\label{z13} \lb \mu =0\ee
 and}
 \bb\label{z14} {\cal L}_{\xi}g_{ij}=\mu g_{ij}.\ee

 \rm 6.2) {\it If $p=2$ and $n=6$, the symmetry group is determined
 by
 \bb\label{z15} X=\xi^{i}(x)\f{\p}{\p x^{i}}+(-\mu (x)\; u +\frac{1}{2}\lb \mu )\frac{\p}{\p
 u} , \ee
 where $\mu $ is a biharmonic function on $M^n$:
 \bb\label{z16} {\lb }^{\!\!\!2}\mu =0 \ee
 and $\xi =\xi^{i}(x){\f{\p}{\p x^{i}}}$ is a conformal Killing
 vector field such that}
 \bb\label{z17} {\cal L}_{\xi}g_{ij}=\mu g_{ij}.\ee

 \

\rm We would like to point out that parts of our results are close
to those in Ibragimov's book \cite{i}. Nevertheless we hope that
they complement them and give useful insights. In fact, the
results on conformally invariant equations established in \cite{i}
are our main motivation to write the present work.

It is clear that the projection of the Lie point symmetries listed
in Theorem 1 on the space of independent variables are, in fact,
special conformal Killing vector fields generating some subgroups
of the conformal group of $(M^{n},g)$, which we shall call special
conformal groups generated by symmetries. Thus the Theorem 1
covers the Ibragimov's theorem \cite{ibr1,i} establishing this
fact.

\

 {\bf Corollary 1. } {\it Let $(M^{n},g)$ be a compact manifold
 without boundary. Then the Lie point symmetry group of
 $(\ref{1.1})$ with an arbitrary $f(u)$ coincides with the
 isometry group $Isom(M^n,g)$.

 If $f(u)=0$ the symmetry group is generated by
 \[X= \xi + (cu+b)\frac{\p}{\p u}, \]
 where $\xi $ is a Killing vector field $($that is ${\cal
 L}_{\xi}g_{ij}=0$$)$, $c$ and $b-$arbitrary constants.

 If $f(u)=u$, then the symmetry generator has the form
 \[ X=\xi + (cu+b(x)) \frac{\p}{\p u},\]
 where $\xi $ is a Killing vector field, $\int_M b(x)dV=0$, $\lb b
 +b=0$ and $c$ is an arbitrary constant.}

 \

\rm By the results in \cite{ya} one can easily obtain the
following estimates on the dimension of the symmetry Lie algebras:

\

{\bf Corollary 2.} {\it Let $\mathfrak{S}$ be a Lie algebra
generated by the symmetries of the nonlinear Poisson equation on
$(M^{n},g)$. Then, the dimension of $\mathfrak{S}$ with an
arbitrary $f(u)$ does not exceed $\f{n(n+1)}{2}$ and the equality
holds if and only if the sectional curvature of $(M^{n},g)$ is
constant.

For some special choices of the function $f(u)$ the dimension of
$\mathfrak{S}$ can be enlarged.}

\begin{enumerate}
\rm\item{\it If $f(u)=k,\,k=const$ or $f(u)=u$, then
$dim(\mathfrak{S})=\infty$ and all finite dimensional subalbegras
possess dimension less than $\f{(n+1)(n+2)}{2}+1$ and the equality
holds if and only if $(M^{n},g)$ is a flat manifold.} \rm\item\it
If $f(u)=e^u$ or $f(u)=u^{p}$, with $p(p-1)(p-\f{n+2}{n-2})\neq
0$, then $dim(\mathfrak{S})\leq\f{n(n+1)}{2}+1$ and the equality
holds if and only if $(M^{n},g)$ is Euclidean. In particular, if
$\mathfrak{p}$ and $\mathfrak{e}$ denotes the symmetry Lie
algebras corresponding to the cases of power and exponential
nonlinearity, then $\mathfrak{p}\approx\mathfrak{e}$, for any
manifold $(M^{n},g)$. \rm\item\it If $f(u)=u^{\f{n+2}{n-2}}$, then
$dim(\mathfrak{S})\leq\f{(n+1)(n+2)}{2}$ and the equality holds if
and only if $(M^{n},g)$ is a flat manifold.
\end{enumerate}\rm

In order not to loose the generality we have not made specific
assumptions on the manifold $M^n$ except $n\geq 3$. (The case
$n=2$ will be treated elsewhere. Some partial results can be found
in \cite{yii}.) Rather we provide a scheme which can be followed
and specialized for any concrete manifold, for which one should
extract further information using its geometrical properties.

Another related point to be emphasized concerns the integrability
conditions for ${\cal L}_{\xi}g_{ij}=\mu g_{ij} $ with $\mu =0$
(isometry), $\mu =const. $ (homothety) or $\mu = \mu (x)$
corresponding to a general conformal transformation (`conformal
motion'). These conditions are in terms of the Riemannian
curvature tensor and have been thoroughly studied in 1950s-60s.
See for instance \cite{ya, sch} and the references therein.
Although we shall not explicitly state the integrability
conditions corresponding to the cases of Theorem 1, we shall
suppose that such conditions hold. (Otherwise the symmetry
determining equations might define the set $\{0\}$. The latter, of
course, may occur for certain manifolds. This simply means that
there are no nontrivial Lie point symmetries of the Poisson
equation (\ref{1.1}) on such manifolds.) Moreover, the number of
the integrability conditions is undetermined for a generic
manifold. This would create another considerable difficulty in
treating the group classification problem in such a general
setting. For this reason we have presented in Theorem 1 just the
relations determining the symmetry groups as special conformal
groups, without entering in differential-geometrical details
regarding $M^n$ like, e.g., positivity, negativity, vanishing or
boundedness of its scalar, sectional or Ricci curvatures as well
as of the respective consequences. For a variety of such results
see \cite{yau1} and the references therein.

Our next purpose in this paper is to find out which of the above
symmetries are variational or divergence
 symmetries.

 \

 {\bf Theorem 2.} \rm 1.) {\it For an arbitrary $f(u)$
 any symmetry of $(\ref{1.1})$ is a variational symmetry, that is,
 the isometry group of $(M^{n},g)$ and the variational symmetry group
 of $(\ref{1.1})$ coincide.}

 \rm 2.) {\it In the exponential case $f(u)=e^u$ the only variational
 symmetries are the isometries of $(M^n,g)$.}

 \rm 3.) {\it In the power case $f(u)=u^p$, $p\neq 0$, $p\neq 1$,
 the symmetry $(\ref{z10})$ is variational if and only if
 \bb\label{z18} p=\frac{n+2}{n-2},\ee
 that is, $p+1$ equals to the critical Sobolev exponent.}

 \rm 3.1) {\it If $p=\frac{n+2}{n-2} $, $n\neq 6$, then the symmetry
 generated by $(\ref{z12})$ is a divergence symmetry.}

 \rm 3.2) {\it If $p=2$ and $n=6$, then the symmetry generated by $(\ref{z15})$ is a
 divergence symmetry.}

 \rm 4.) {\it In the linear cases we have:}

 \rm 4.1) {\it If $f(u)=0$, the symmetry $(\ref{z2})$ is a Noether symmetry if and
 only if $c=0$.}

\rm 4.2) {\it If $f(u)=k$, $k\neq 0$, the symmetry $(\ref{h1})$ is
a Noether symmetry if and only if $b=c=0$.}

\rm 4.3) {\it If $f(u)=u$, the symmetry $(\ref{z2})$ is a Noether
symmetry if and only if $c=0$.} \rm

 \

\rm From case 4.2) we conclude that the Noether symmetry group for
the non-homogeneous case $f(u)=k$ coincides with the isometry
group of $(M^{n},g)$. Also from case 4.3), if $c=0$ in (\ref{z2}),
then, except the term $b\f{\p}{\p u}$ in (\ref{z2}), the Noether
symmetry group of the homogenous case coincides with the symmetry
group of the critical case (\ref{z12}).

We would also like to observe that in the critical cases 3.1) and
3.2) of Theorem 2 all Lie point symmetries are Noether symmetries.
The general property stating that a Lie point symmetry of an
equation (or system) is a Noether symmetry if and only if the
equation parameters assume critical values has been established
and discussed in \cite{yb1}. Recall that, as it is well known, the
so-called critical exponent is found as the critical power for
embedding theorems of Sobolev type. It is also related to some
numbers dividing the existence and nonexistence cases for the
solutions of differential equations, in particular semilinear
differential equations with power nonlinearities involving the
Laplace operator. The above mentioned property traces a connection
between these two notions: the Noether symmetries and the
`criticality' of the equation. It relates two important theorems,
namely, the Sobolev Theorem and the Noether Theorem. Theorem 2
shows that this property holds also in the context of Riemannian
manifolds. In fact, this is another motivation to write the
present paper. In this regard, we can see that the widest symmetry
group admitted by the Poisson equation may be the full conformal
group of $M^n$. Namely:

\

{\bf Corollary 3.} {\it If the scalar curvature of $(M^n,g)$,
$n\geq 3$, vanishes, then the widest symmetry group admitted by
the Poisson equation $(\ref{1.1})$ is achieved for the critical
equation
\[ \lb u+u^{(n+2)/(n-2)}=0.\]
In this case it coincides with the conformal group of} $(M^n,g)$.

\

\rm The latter equation and its invariant properties have been
studied in \cite{i}.

We remark that if the manifold is flat, then the symmetry group of
the nonlinear cases is maximal in the critical case.

Now we shall state the conservation laws corresponding to the
found Noether symmetries. Before doing this it is worth mentioning
that there are powerful modern methods to obtain conservation laws
due to George Bluman et al.
\cite{ab1,ab2,ab3,a,i1,ibjmaa,k1,naz,wolf}. We believe that these
methods can be very useful in the study of various differential
geometric problems. However, we have chosen here the classical
approach since we have at our disposal explicit formulae for the
potentials ensured by the Noether's theorem (see sections 8 and
9), whose determination is usually the major difficulty. Thus it
is immediate, simple and natural to apply the classical approach.

\

{\bf Corollary 4.} {\it The conservation laws corresponding to the
Noether symmetries of equation $(\ref{1.1})$ are classified as
follows:}

\begin{enumerate}
\rm\item {\it If $f(u)$ is an arbitrary function, then the
conservation law is $D_{i}A^{i}=0,$ where \bb\label{cc1}
A^{k}=\sqrt{g}(\f{1}{2}g^{ij}\xi^{k}-g^{kj}\xi^{i})u_{i}u_{j}-\sqrt{g}\xi^{k}F(u),\ee
and $X=\xi^{i}\f{\p}{\p x^{i}}$ is a Killing vector field on
$(M^{n},g)$.}

\rm\item {\it If $f(u)=0$, then the conservation law is
$D_{i}A^{i}=0$, where \bb\label{cc4} A^{k}  =
\sqrt{g}(\f{1}{2}g^{ij}\xi^{k}-g^{kj}\xi^{i})u_{i}u_{j}
+\f{2-n}{4}\sqrt{g}g^{kj}(\mu u
u_{j}-\f{1}{2}\mu_{j}u^{2})+\sqrt{g}g^{jk}(bu_{j}-b_{j}u), \ee and
$\xi $, $b$ $\mu $ satisfy} $(\ref{z3}), (\ref{z4}), (\ref{z5})$.

\rm\item {\it If $f(u)=k,\,k=const.$, then the conservation law is
given by equation $(\ref{cc4})$ with }$b=0$.

\rm\item {\it If $f(u)=u$, then the conservation law is
$D_{i}A^{i}=0$, where \bb\label{cc5} \ba{l c l}
A^{k}  & = & \ds{ \sqrt{g}(\f{1}{2}g^{ij}\xi^{k}-g^{kj}\xi^{i})u_{i}u_{j} +\f{2-n}{4}\sqrt{g}g^{kj}(\mu u u_{j}-\f{1}{2}\mu_{j}u^{2})}\\
\\
&
&\ds{+\sqrt{g}g^{jk}(bu_{j}-b_{j}u)-\f{1}{2}\xi^{k}\sqrt{g}u^{2}}
\ea \ee and $\xi $, $b$ $\mu $ satisfy} $(\ref{z5}), (\ref{z6}),
(\ref{z7})$.

\rm\item {\it If $n\neq 6$ and $f(u)=u^{\f{n+2}{n-2}}$, then the
conservation law is $D_{i}A^{i}=0,$ where \bb\label{cc2} \ba{l c
l} \ds{A^{k}}&  =&
\ds{\sqrt{g}(\f{1}{2}g^{ij}\xi^{k}-g^{kj}\xi^{i})u_{i}u_{j}+\f{2-n}{2n}\sqrt{g}\xi^{k}u^{\f{2n}{n-2}}}\\
\\
&&\ds{+\f{2-n}{4}\sqrt{g}g^{kj}(\mu u
u_{j}-\f{1}{2}\mu_{j}u^{2})}, \ea \ee and $\xi $, $\mu $ satisfy}
$(\ref{z13}), (\ref{z14})$.

\rm\item {\it If $n=6$ and $f(u)=u^{2}$, then the conservation law
is $D_{i}A^{i}=0$, where \bb\label{cc3} \ba{l c l}
\ds{A^{k}} & = & \ds{\sqrt{g}(\f{1}{2}g^{ij}\xi^{k}-g^{kj}\xi^{i})u_{i}u_{j}-\f{1}{3}\sqrt{g}\xi^{k}u^{3}}\\
\\
 & & \ds{ +\sqrt{g}g^{kj}[\f{1}{2}(\lb\mu u_{j}+\mu_{j}u^{2})-(\mu u u_{j}+(\lb\mu)_{j}u)]}
\ea \ee and $\xi $, $\mu $ satisfy} $(\ref{z16}), (\ref{z17})$.

\end{enumerate}

\rm For some applications of symmetries and conservation laws see
\cite{bk,yb2,yi2,i,igor}.

This paper is organized as follows. Section 2 includes the
geometric preliminaries and introduces notations and conventions
used in this paper. Further, the determining equations for the
symmetry coefficients are obtained in the section 3. The
connections between isometry groups and symmetry groups are
established in the section 4. The group classification for the
linear cases is obtained in the section 5 and for the nonlinear
cases in sections 6 and 7. The proof of Corollary 1 concerning the
Lie point symmetries in the case of compact manifolds without
boundary is presented in the section 8. The Noether symmetries are
found in sections 9 and 10. In order to illustrate the main
results, some examples are presented in the section 11 in which we
perform the group classification and establish the Noether
symmetries and their respective conservation laws for the
nonlinear Poisson equations in the Thurston geometries, namely:
$\mathbb{R}^{3}$, the three-dimensional Hyperbolic space $
\mathbb{H}^3$, the sphere $\mathbb{S}^{3}$, the three-dimensional
solvable group Sol, the product spaces $\mathbb{S}^2\times\R$,
$\mathbb{H}^2\times\R$, the universal covering of ${SL_2}(\R)$ and
the three-dimensional Heisenberg Group.

\

\section{Preliminaries}

\

In this section we introduce notation. We also state some results
which will be used later.

The Riemann tensor of $(M^{n},g)$ is given by \bb\label{2.1}
R^{i}_{\;jks}=\G^{i}_{jk,s}-\G^{i}_{js,k}+\G^{i}_{ls}\G^{l}_{jk}-\G^{i}_{lk}\G^{l}_{js}.
\ee

The Ricci tensor \bb\label{2.2} R^{i}_{\;s}=g^{jk}R^{i}_{\;jks}
\ee and its trace $R:=R^{i}_{\;i}$ is the scalar curvature of $M$.

For any contravariant vector field $T=(T^i)$ the following
commutation relation holds \bb\label{2.3}
(\nabla_{k}\nabla_{l}-\nabla_{l}\nabla_{k})T^{i}=-R^{i}_{\;
skl}T^{s}. \ee See \cite{dnf}.

We observe that the Riemann and Ricci tensors used in this paper
coincide with those in Yano's book \cite{ya} and in Dubrovin,
Fomenko and Novikov's book \cite{dnf}; they are negatives of the
respective tensors used by Ibragimov in \cite{i}.

We shall need some auxiliar results which are presented in a
sequence of lemmas.

\

{\bf Lemma 1.} {\it If $\xi$ is a conformal Killing vector field
satisfying \bb\label{2.4} \nabla^{i}\xi^{j}+\nabla^{j}\xi^{i}=\mu
g^{ij}, \ee then the covariant divergence} \bb\label{2.6} div
(\xi)=\nabla_{j}\xi^{j}=\f{n}{2}\mu. \ee

\

{\bf Proof.} \rm Just take the trace in $(\ref{2.4})$.

\

{\bf Lemma 2.} {\it If $\xi$ is a conformal Killing vector field
$($see equation $(\ref{2.4})$$)$ then} \bb\label{2.5} \lb
\xi^{i}+R^{i}_{\;j}\xi^{j}=\f{2-n}{2}g^{ij}\mu_{j}. \ee

\

{\bf Proof.} \rm Applying the covariant derivative operator
$\nb_j$ to equality $(\ref{2.4})$ and summing up we obtain
 \bb \label{a1} \nb_j \nabla^{i}\xi^{j}+\lb \xi^{i}=
 g^{ij} {\mu }_j.
 \ee
 On the other hand, from $(\ref{2.3})$ it follows that
 \[\nb_j \nabla^{i}\xi^{j}- \nb^i\nb_j \xi^{j} =R^{i}_{\;s}\xi^{s}  \]
 and hence
 \bb \label{a2} \nb_j \nabla^{i}\xi^{j} = \nb^i div (\xi ) +
 R^{i}_{\;s}\xi^{s} . \ee
 Then $(\ref{2.6})$, $(\ref{a1})$ and $(\ref{a2})$ imply the
 relation $(\ref{2.5})$.

 \

 {\bf Lemma 3.} \rm (Yano \cite{ya}). {\it If $\xi$ is a conformal vector
 field such that \bb\label{2.7} {\cal L}_{\xi}g_{ij}=\mu g_{ij},
 \ee then} \bb\label{2.8} \lb \mu=-\f{1}{n-1}(\L_{\xi}R+\mu R). \ee

\

\rm The next two lemmas concern the form of the semilinear Poisson
equation (\ref{1.1}) and its variational structure.

\

{\bf Lemma 4.} {\it The Poisson equation $(\ref{1.1})$ can be
written in the equivalent form} $(\ref{a3})$.

\

{\bf Proof.} \rm This can be seen by performing explicitly the
partial differentiations in the Laplace-Beltrami operator and
using the formula \bb\label{a4} (\sqrt{g} g^{ik})_{,k}= -
g^{pq}\G_{pq}^{i} \sqrt{g}.\ee

 \

 {\bf Lemma 5.} {\it The Poisson equation $(\ref{1.1})$ has a variational structure and
 it is (formally)
 the Euler-Lagrange equation of a functional $\int_{M} L dx $, where
 the Lagrangian
 \bb\label{2.9} L=\f{\sqrt{g}}{2}g^{ij}u_{i}u_{j}-F(u)\sqrt{g}, \ee
 and} $F'(u)=f(u)$.

 \

 {\bf Proof.} \rm First we apply to $L$ the Euler operator
 \[ E=\frac{\kd }{\kd u} - D_k\frac{\kd }{\kd u_k}, \]
 where
 \[D_k=\frac{\kd  }{\kd x^k}+ u_k\frac{\kd }{\kd u} + u_{k s}\frac{\kd }{\kd u_s}+...\]
 is the total derivative operator. Then, after simplifying, the
 equation $(\ref{1.1})$ is obtained.

 \

 \section{The determining equations}

\

\rm In this section we shall obtain the set of linear partial
differential equations determining the Lie point symmetries of the
Poisson equation (\ref{1.1}).

 To begin with, let
 \[ X={\xi }^i (x,u) \frac{\kd }{\kd x^i}+ {\eta } (x,u) \frac{\kd }{\kd u } \]
 be a partial differential operator on $M^n\times {\mathbb{R}}$
 which is infinitesimal generator of such a symmetry. Let
 $X^{(1)}$ be the first order prolongation of $X$. (See \cite{ba, bk, i, ol, ov} for the
 corresponding definitions.) First we shall simplify the form of
 $X$.

 \

{\bf Proposition 1.} {\it If $n\geq 2$, then the infinitesimals of
the symmetry $X$ take the form }\bb\label{3.1} \left\{ \ba{l c l}
\xi^{i} & = & \xi^{i}(x),\\
\\
\eta & = & a(x)\;u+b(x), \ea\right. \ee {\it where $a=a(x)$ and
$b=b(x)$ are smooth functions.}

\

{\bf Proof.} \rm This proposition follows from two theorems of
Bluman \cite{gb,bk}.

\

\rm The following intermediate result, up to notation, is the same
as in Ibragimov's book \cite{i}, pp. 114-115.

\

{\bf Proposition 2.} \rm (Ibragimov \cite{i}). {\it The
infinitesimals of the Lie point symmetries of equation
$(\ref{1.1})$ satisfy the relations:} \bb\label{3.2}
\xi^{k}g^{ij}_{\;\;,k}-g^{ik}\xi^{j}_{\;,k}-g^{jk}\xi^{i}_{\;,k}+a
g^{ij}=\lambda g^{ij}, \ee

\bb\label{3.3} 2g^{ij}
a_{j}-g^{jk}\xi^{i}_{\;,jk}+\G^{j}\xi^{i}_{\;,j}-\G^{i}_{\;,j}\xi^{j}-a
\G^{i}=-\lambda \G^{i}, \ee

\bb\label{3.4} (\lb a)\;u+(a \;u+b)\;f'(u)+\lb b = \lambda f(u),
\ee where $\lambda=\lambda(x)$ and $\G^{i}:=g^{pq}\G^{i}_{pq}$.

\

{\bf Proof.} \rm This proposition follows from Proposition 1 and
the definition of Lie point symmetry.

\

{\bf Proposition 3.} {\it Let $\xi=\xi^{i}\f{\p}{\p x^{i}}$. Then
the determining equations can be written in the following
equivalent form:}

\bb\label{3.5}
\nabla_{i}\xi_{j}+\nabla_{j}\xi_{i}=\L_{\xi}g_{ij}=\mu g_{ij}, \ee

\bb\label{3.6} 2g^{ij}a_{j}-(\lb \xi^{i}+R^{i}_{\;j}\xi^{j})=0,
\ee

\bb\label{3.7} a uf'(u)+b f'(u)-\lambda f(u)+(\lb a) u +\lb b=0,
\ee {\it where $\mu:=a-\lambda$ and $R^{i}_{\;j}$ is the Ricci
tensor of }$g$.

\

{\bf Corollary 5.} {\it  The relation $(\ref{3.6})$ is equivalent
to \bb\label{3.8} a_{i}=\f{2-n}{4}\mu_{i} \ee which, itself, is
equivalent to \bb\label{3.9} \lambda_{i}=\f{n+2}{n-2}a_{i}. \ee
Thus, the determining equations are $(\ref{3.5})$, $(\ref{3.8})$
and} \bb\label{h3} a\; uf'(u)+b\; f'(u)+(\mu
-a)\;f(u)+\frac{2-n}{4}(\lb \mu )\; u +\lb b=0.\ee

\

{\bf Corollary 6.} {\it The determining equations are:}
\bb\label{b3}
\nabla_{i}\xi_{j}+\nabla_{j}\xi_{i}=\L_{\xi}g_{ij}=\mu g_{ij}, \ee
\bb\label{b4} a_{i}=\f{2-n}{4}\mu_{i} \ee \bb\label{b5} a\;
uf'(u)+b\; f'(u)+(\mu -a)\;f(u)+\frac{n-2}{4(n-1)}\;[ {\xi}^i
R_{,i}+\mu R]\;u +\lb b=0. \ee

 \

 \rm The statement of Corollary 6 is explicitly announced in \cite{i},
 pp. 115 - 116, without a detailed proof. For the sake of clearness and completeness we have
 decided to present here the corresponding proof, dividing the
 procedure in three steps: Proposition 3, Corollary 5 and Corollary
 6.

 \

 {\bf Proof of Proposition 3.} \rm The equivalence between
(\ref{3.4}) and (\ref{3.7}) is obvious. The equivalence between
(\ref{3.2}) and (\ref{3.5}) is clear from the definitions of Lie
derivative and conformal Killing vectors. Indeed, from
(\ref{3.2}):
$$
\xi^{k}g^{ij}_{\;\;\;,k}-g^{is}\xi^{j}_{\;,s}-g^{js}\xi^{i}_{\;,s}=-\mu
g^{ij}.
$$
The formula
$$ g^{ij}_{\;\;\;,k}=-\G_{ks}^{i}g^{sj}-\G^{j}_{ks}g^{is},$$
substituted above, implies
$$\nabla ^{i}\xi^{j}+\nabla^{j}\xi^{i}=\mu g^{ij}.$$
Hence and by the definition of Lie derivative, it follows that
(\ref{3.5}) holds. And vice versa: (\ref{3.5}) implies
(\ref{3.2}).

Further, we shall prove the equivalence of (\ref{3.3}) and
(\ref{3.6}). From (\ref{3.3}) we have \bb\label{4.1}
2g^{ij}a_{j}-(g^{jk}\xi^{i}_{\;\;,jk}-\G^{j}\xi^{i}_{,j}+\G^{i}_{,j}\xi^{j}+\mu\G^{i})=0.
\ee

We denote
$$A^{i}:=g^{jk}\xi^{i}_{\;\;,jk},\;\;\; B^{i}:=\G^{j}\xi^{i}_{\;,j},\;\;\;
C^{i}:=\G^{i}_{\;,j}\xi^{j}.$$ Then (\ref{4.1}) reads
\[ 2g^{ij}a_{j}-(A^{i}-B^{i}+C^{i}+\mu\G^{i})=0.\]
We aim to express $A^{i},\; B^{i}$ and $C^{i}$ in the terms of
covariant derivatives. From \bb\label{4.3}
\nabla_{j}\xi^{i}=\xi^{i}_{\;,j}+\G^{i}_{js}\xi^{s} \ee and
\bb\label{4.2} \nabla_{j}\G^{i}=\G^{i}_{\;,j}+\G_{js}^{i}\G^{s}
\ee we obtain \bb\label{4.4}
B^{i}=\G^{j}(\nabla_{j}\xi^{i}-\G^{i}_{js}\xi^{s})=\G^{j}\nabla_{j}\xi^{i}-\G^{j}\G_{js}^{i}\xi^{s}
\ee and \bb\label{4.5}
C^{i}=(\nabla_{j}\G^{i}-\G^{i}_{js}\G^{s})\xi^{j}=(\nabla_{j}\G^{i})\xi^{j}-\G^{i}_{js}\G^{s}\xi^{j}.
\ee Further:
$$\nabla_{k}\xi^{i}_{\;,j}=\xi^{i}_{\;,jk}+\G^{i}_{kl}\xi^{l}_{\;,j}
-\G^{l}_{kj}\xi^{i}_{\;,l}$$

\bb\label{4.6} \ba{l c l} \xi^{i}_{\;,jk}& = &
\nabla_{k}(\nabla_{j}\xi^{i}-\G^{i}_{js}\xi^{s})-
\G^{i}_{kl}(\nabla_{j}\xi^{l}-\G^{l}_{js}\xi^{s})+
\G^{l}_{kj}(\nabla_{l}\xi^{i}-\G^{i}_{ls}\xi^{s})\\
\\
& = & \nabla_{k}\nabla_{j}\xi^{i}-(\nabla_{k}\G^{i}_{js})\xi^{s}-
\G^{i}_{js}\nabla_{k}\xi^{s}-\G^{i}_{kl}\nabla_{j}\xi^{l}\\
\\
& + &
\G^{i}_{kl}\G^{l}_{js}\xi^{s}+\G^{l}_{kj}\nabla_{l}\xi^{i}-\G^{i}_{ls}\xi^{s}\G^{l}_{kj}.
\ea \ee

\

From (\ref{4.1}), (\ref{4.4}), (\ref{4.5}) and (\ref{4.6}) it
follows that \bb\label{4.7} \ba{c}
2g^{ij}a_{j}-\lb\xi^{i}+g^{jk}(\nabla_{k}\G^{i}_{js})\xi^{s}+g^{jk}\G^{i}_{js}\nabla_{k}\xi^{s}+\G^{i}_{kl}g^{jk}\nabla_{j}\xi^{l}-g^{jk}\G^{i}_{kl}\G^{l}_{js}\xi^{s}\\
\\
-\G^{l}\nabla_{l}\xi^{i}+\G^{l}\G^{i}_{ls}\xi^{s}-(\nabla_{j}\G^{i})\xi^{j}+\G^{s}\G^{i}_{js}\xi^{j}+\G^{j}\nabla_{j}\xi^{i}-\G^{j}\G^{i}_{js}\xi^{s}-\mu\G^{i}=0.
\ea \ee On the other hand, from (\ref{3.5})
$$\nabla^{k}\xi^{l}=-\nabla^{l}\xi^{k}+\mu g^{lk}$$
and hence \bb\label{4.8}
-\G^{i}_{kl}g^{jk}\nabla_{j}\xi^{l}=\G^{i}_{kl}\nabla^{l}\xi^{k}-\mu
g^{lk}\G^{i}_{lk}. \ee Substituting (\ref{4.8}) into (\ref{4.7})
we obtain after renaming the indices and cancelling some terms:
\bb\label{4.9} 2g^{ij}a_{j}-(\lb \xi^{i}+w^{i}_{s}\xi^{s})=0, \ee
where \bb\label{4.10}
w^{i}_{s}=-g^{jk}(\nabla_{k}\G^{i}_{js})+g^{jk}\G^{i}_{kl}\G^{l}_{js}+\nabla_{s}\G^{i}-\G^{l}\G^{i}_{sl}.
\ee Then we expresses the covariant derivatives in (\ref{4.10}) in
the terms of usual partial derivatives. In this way we get that
\bb\label{4.11}
w^{i}_{s}=g^{jk}(\G^{i}_{jk,s}-\G^{i}_{js,k}+\G^{i}_{ls}\G^{l}_{jk}-\G^{i}_{lk}\G^{l}_{js})=
g^{jk}R^{i}_{\;jks}=R^{i}_{\;s}, \ee where $R^{i}_{\;jks}$ and
$R^{i}_{\;s}$ are the components of the Riemann and Ricci tensors
respectively. Thus (\ref{4.9}), (\ref{4.10}) and (\ref{4.11})
imply (\ref{3.6}).

\

{\bf Proof of Corollary 5.} \rm From (\ref{2.5}) and (\ref{3.6})
we obtain
$$2g^{ij}a_{j}-\f{2-n}{2}g^{ij}\mu_{j}=0,$$
that is,
$$a_{i}=\f{2-n}{4}\mu_{i}.$$
The rest of the proof is straightforward.

\

{\bf Proof of Corollary 6.} \rm The conclusion follows from
Corollary 4 and Lemma 3 (see (\ref{2.8})).

\

\section{The isometry group and the Lie point symmetry group for arbitrary $f(u)$}

\

\rm In this section we prove the first part of Theorem 1.

Let $X$ be a symmetry of (\ref{1.1}).  Then $X$ has necessarily
the form given in Proposition 1. From (\ref{3.7}), equating to
zero the terms involving $u$, we obtain:
$$a=b=\lambda=\mu=0.$$
Hence, $\eta=0$ and $X=\xi^{i}\f{\p}{\p x^{i}}$. From (\ref{3.5})
with $\mu=0$, it follows that $X$ is an isometry.

Let $X=\xi^{i}\f{\p}{\p x^{i}}$ be an infinitesimal isometry of
$M^{n}$. Therefore $\L_{X}g_{ij}=0$. Hence the equation
(\ref{3.5}) holds with $\mu=0$. By the form of $X$, $\eta=0$ and
thus $a=b=0$. Hence the relation (\ref{h3}) is satisfied.
Obviously (\ref{3.8}) is satisfied since $a=\mu=0$. Therefore
$X=\xi^{i}\f{\p}{\p x^{i}}$ is a Lie point symmetry.

In this way we have proved that
\[ X={\xi }^i (x,u) \frac{\kd }{\kd x^i}+ {\eta } (x,u) \frac{\kd }{\kd u }  \]
 is a Lie point symmetry of (\ref{1.1}) with arbitrary $f(u)$ if
 and only if $\eta =0$ and $X={\xi }^i (x) \frac{\kd }{\kd
 x^i}=\xi $ is an infinitesimal isometry of $M^n$. In other words,
 the isometry group of $(M^n,g)$ is the invariance group of its
 Laplace-Beltrami operator!

 In this case, by (\ref{2.8}) we have ${\cal L}_{\xi } R=0$, an integrability
 condition for ${\cal L}_{\xi } g_{ij}=0$, which we suppose holds
 true. See \cite{ya} for details concerning the integrability
 conditions for $ {\cal L}_{\xi}g_{ij}=0$ or, more generally,
 ${\cal L}_{\xi}g_{ij}=\mu g_{ij}$.

 \

 \section{The Lie point symmetries in the linear cases}

 \

 In this section we prove parts 2.), 3.) and 4.) of Theorem 1.

 Let $f(u)=0$. Then from (\ref{h3}) we conclude that $\lb b=0$ and
 $\lb \mu =0$. Integrating (\ref{3.8}) we obtain $a=\f{2-n}{4}\mu
 +c$, where $c$ is an arbitrary constant. This completes the proof
 of the second part of Theorem 1.

 Further, let $f(u)=k=const\neq 0$. From (\ref{h3}) we have
 \[ (\mu -a)k + \f{2-n}{4}(\lb \mu )\;u +\lb b=0.\]
 Hence, equating to zero the coefficient of $u$ and the free term,
 we obtain  \bb\label{h4} \lb\mu =0\ee and \bb\label{h5}(\mu -a)k  +\lb b=0.\ee
 Again from (\ref{3.8}) we get $a=\f{2-n}{4}\mu
 +c$, where $c$ is an arbitrary constant. Substituting this
 expression for $a$ in (\ref{h5}) we find
 \[\mu = \f{4}{(n+2)}(c-\f{1}{k}\lb b). \]
 It remains to put the latter into (\ref{h4}) to conclude that ${\lb }^{\!\!\!2} b =0
 $.

 Let now $f(u)=u$. From (\ref{h3}) we have
 \[ a\;u +b + (\mu -a)\;u + \f{2-n}{4}(\lb \mu )\;u +\lb b=0.\]
 Hence, equating to zero the coefficient of $u$ and the free term,
 we obtain  $\lb b+b=0$ and $\f{2-n}{4}\lb \mu +\mu=0$. This
 completes the proof of Theorem 1 in the linear cases.

 \

 \section{The Lie point symmetries in the case of exponential
 nonlinearity}

 \

 Let $f(u)=e^u$. From (\ref{h3}) we obtain
 \[ a\;ue^u +(b+\mu -a)\; e^u +\f{2-n}{4}(\lb \mu )\;u +\lb b=0.\]
 Hence, equating to zero the coefficients of $ue^u$, $e^u$, $u$  and the free term,
 we obtain  $a=0$, $b=-\mu $ and $\lb b=\lb \mu=0$. From (\ref{3.8})
 \[ \f{2-n}{4} {\mu }_i =0\]
 since $a=0$. Thus ${\mu }_i =0$ because $n\geq 3$. Hence $\mu
 =const.$ This completes the proof.

 We observe that this case, in fact, is a Liouville-Gelfand problem
 on Riemannian manifolds.

 \

 \section{The Lie point symmetries in the case of power nonlinearity}

 \

 In this section we shall prove parts 6.), 6.1) and 6.2) of
 Theorem 1.

 Let $f(u)=u^p$. The cases $p=0$ and $p=1$ have already been
 considered and we may suppose that $p\neq 0$ and $p\neq 1$.

 From (\ref{h3}) we obtain
 \[ [(p-1)a +\mu ]\; u^p +pb\;u^{p-1}+\f{2-n}{4}(\lb \mu )\;u +\lb b=0.\]
 Let $p\neq 2$. Then, equating to zero the coefficients of $u^p$, $u^{p-1}$, $u$
 and the free term, implies that
 \bb\label{z22} a=\f{1}{1-p}\mu , \ee
 $pb=0$ and $\lb\mu=0$. Since $p\neq 0$ it follows that $b=0$.
 From (\ref{3.8}) and (\ref{z22}) it follows that
 \bb\label{z23} \left [ \f{1}{1-p}+\f{n-2}{4} \right ]\;{\mu }_i=0 .\ee

 (i) If $p\neq \f{n+2}{n-2}$, then (\ref{z23}) implies ${\mu
 }_i=0$ for all $i$ and thus $\mu=const.$ and $\xi $ is a
 homothety.

 (ii) Let $p = \f{n+2}{n-2}$ and $n\neq 6$. (Otherwise $p=2$.)
 From (\ref{z22}) it is clear that the symmetry has the form
 announced in 6.1) of Theorem 1.

 Let now $p=2$. From (\ref{h3}) we have
 \bb\label{z24} (a+\mu )\; u^2 + (2b+\f{2-n}{4}(\lb \mu ))\;u +\lb
 b=0.\ee
 Hence
 \bb\label{z25} a=-\mu ,\ee
 \bb\label{z26} 2b+\f{2-n}{4}(\lb \mu )=0 \ee
 and
 \bb\label{z27} \lb b=0. \ee
 From (\ref{3.8}) and (\ref{z25}) it follows that
 \bb\label{z28} \f{6-n}{4}{\mu }_i=0. \ee

 Again we have to consider two cases.

 I.) If $n\neq 6$ then from (\ref{z28}) we obtain $ {\mu }_i=0$
 for all $i$ and thus $\mu=const.$ Hence and from (\ref{z26}) it
 follows that $b=0$. We see that this case is included in case 6.)
 of Theorem 1.

 II.) Let $n=6$. Now we cannot conclude from (\ref{z28}) that $\mu
 $ is a constant. From (\ref{z26}) we express $b$ as a function of
 $\mu $:
 \[ b= \f{n-2}{8}\lb\mu = \f{1}{2}\lb\mu .\]
 Substituting this expression for $b$ into (\ref{z27}) we obtain ${\lb }^{\!\!\!2}\mu =0
 $. This completes the proof of Theorem 1.

 \

 {\bf Proof of Corollary 3.} \rm If $R=0$ then by Lemma 3,
 (\ref{2.8}) any conformal transformation (\ref{2.7}) satisfies
 $\lb \mu =0$ and hence the conclusion follows immediately.

 \

 \

\section{The Lie point symmetries in the case of compact manifolds without
boundary}

\

In this section we prove Corollary 1.

1.) Let $f(u)$ be an arbitrary function. By Theorem 1 the symmetry
group coincides with the isometry group of $(M^n,g)$.

2.) If $f(u)=0$ by Theorem 1 the symmetry is determined by (4),
(5), (6), (7). Since $\mu $ and $b$ are harmonic, by the E. Hopf's
maximum principle it follows that $\mu = const.$ and $b=const.$
But $\mu = 2 div(\xi)/n $ (see Lemma 1). Hence
\[\mu Vol(M^n)= \int_{M^n} \mu \;dV = \frac{2}{n} \int_{M^n} div(\xi) \;dV = 0\]
by the Green's theorem. Thus $\mu =0$.

3.) Here we observe that the Poisson equation with $f(u)=k$ on
compact manifolds without boundary makes sense only if $k=0$. This
follows directly from the Green's theorem.

The proof of Corollary 1 in the rest of the cases in Theorem 1 is
similar to the above presented. For this reason we shall not
present further details merely pointing out that the constancy of
the conformal factor $\mu $, the harmonicity of $\mu $ or
biharmonicity of $\mu $ imply that $\mu=0$ by the maximum
principle and the Green's theorem.

 \section{The Noether symmetries in nonlinear cases}

 \

 In this and the next sections we present the proof of Theorem 2 divided in
 several propositions and lemmas.

 In order to apply the infinitesimal criterion for invariance \cite{ol}, p.
 257, we need the following

 \

 {\bf Proposition 4.} {\it Let
 \bb\label{n1} X=\xi^{i}(x)\f{\p}{\p x^{i}} +[a(x)\; u +b(x)]\frac{\p}{\p
 u}, \ee
 where $a$, $b$ and $\xi^{i}$ are smooth functions, be a partial differential
 operator on $M^n\times \mathbb{R}$. Then}

 \bb\label{n2} \begin{array}{lll} X^{(1)}L + L\;D_i{\xi }^i& = & \f{1}{2}[ g^{ks}div(\xi )+2 a g^{ks}
 -\nabla^{k}\xi^{s}-\nabla^{s}\xi^{k}]
 \sqrt{g} u_k u_s \\
 & & \\
 &- &  \sqrt{g}\;div(\xi )\; F(u)-\sqrt{g}a\;uf(u) -\sqrt{g}b\;f(u)\\
 \\
 & +& (a_i\;u+b_i)\sqrt{g}g^{is}u_s, \end{array} \ee
 {\it where $X^{(1)}$ is the first order prolongation of $X$ and
 the function $L$ is given in} $(\ref{2.9})$.

 \

 {\bf Proof.} \rm By a straightforward calculation we obtain the
 first order prolongation
 \[X^{(1)} = X + \left ( a_i u +b_i +(a {\de }^j_i- {{\xi }^j}_{,i})u_j \right )
 \frac{\kd }{\kd u_i } \]
 where `,' means partial derivative: ${{\xi }^j}_{,i}=\kd {\xi }^j /\kd x^i
 $ and ${\de }^j_i $ is the Kronecker symbol. We apply
 $X^{(1)}$ to $L$ given by (\ref{2.9}) and then change some of the indices in the
 obtained expression. In this way we get that
 \bb\label{t1} \begin{array}{lll}
 X^{(1)} L + L\ds{\frac{\kd {\xi }^i}{\kd x^i} }& = &
  \ds{ \left [ \frac{1}{2} {\xi
 }^i (g^{ks} \sqrt{g})_{,i} + \left ( a {\de }^k_i-
 {{\xi }^k}_{,i} \right )
 g^{is} \sqrt{g} \right. }+
 \ds{ \left. \frac{1}{2}{{\xi }^i}_{,i} g^{ks} \sqrt{g}
 \right ] u_k u_s }\\
 & & \\
 & -& a u f(u)\sqrt{g}-b f(u)\sqrt{g} + (a_iu+b_i)u_s g^{is}\sqrt{g}\\
 & & \\
 &-& F(u){\xi }^i (\sqrt{g})_i-F(u){{\xi }^i}_{,i}\sqrt{g} . \end{array}
 \ee
 Further we shall make use of the formulae
 \bb\label{h6} (g^{ks} \sqrt{g})_{,i} = - g^{sl} {\Ga }^k_{li} - g^{kl} {\Ga }^s_{li},
 \;\;\;\;\;\;
 (\sqrt{g})_{,i} = {\Ga }^k_{ik} \sqrt{g} ,\ee
 where ${\Ga }$'s are the Christoffel  symbols. Then by the definition
 of the covariant derivative operators ${\nabla }^i$,
 corresponding to the Levi-Civita connection $\nabla $, and the second formula in (\ref{h6}),
 the last two terms in (\ref{t1}) can be written as
 \bb\label{h44} - {\xi }^i
 (\sqrt{g})_{,i} F(u) - {{\xi }^i}_{,i} \sqrt{g} F(u) =-div(\xi ) F(u) \sqrt{g}
 \ee
(We recall that $div(\xi )={\nabla }_i {\xi }^i $ is the covariant
 divergence of $\xi $.)

 Now we denote by $A$ the expression in the right-hand side of the
 first line of
 (\ref{t1}) containing $ u_k u_s $. Using (\ref{h6}) we obtain
 that
 \bb\label{h7} \begin{array}{lll} A &= & {\left \{ \frac{1}{2}{\xi }^i \sqrt{g}
 \left [ \ds{- g^{sl} {\Ga }^k_{li} - g^{kl} {\Ga }^s_{li} + g^{ks}{\Ga }^l_{il}} \} \right ]+
 a g^{ks} \sqrt{g} \right. }\\
& & \\
 &-& \left.  \frac{1}{2}g^{is} \sqrt{g} {{\xi }^k}_{,i}
 - \frac{1}{2}g^{ki} \sqrt{g}{{\xi }^s}_{,i}+\frac{1}{2}g^{ks} \sqrt{g}
 {{\xi }^i}_{,i} \right \} u_k u_s \\
& & \\
 &=& \left \{ - \frac{1}{2}(g^{is} {{\xi}^k}_{,i} + g^{sl} {\Ga }^k_{li} {\xi }^i) -
 \frac{1}{2}( g^{ki} {{\xi }^s}_{,i}+g^{kl} {\Ga }^s_{li} {\xi }^i ) \right. \\
 & & \\
 & +& \left. \frac{1}{2}g^{ks} ( {{\xi }^i}_{,i}+{\Ga }^l_{il} {\xi }^i)+
 a g^{ks} \right \} \sqrt{g}u_k u_s .
 \end{array}\ee
 From the definition
 of the covariant derivative we have that
 \[\begin{array}{lll} {\nb }^s {\xi }^k & =& g^{is} {{\xi}^k}_{,i} + g^{sl} {\Ga }^k_{li}
{\xi }^i, \\
& & \\
 {\nb }^k {\xi }^s &=&g^{ki} {{\xi }^s}_{,i}+g^{kl}
{\Ga }^s_{li} {\xi }^i,\\ & & \\{\nb }_i  {\xi }^i &=& {{\xi
}^i}_{,i}+{\Ga }^l_{il} {\xi }^i.\end{array} \] We substitute
these formulae into (\ref{h7}). Thus
 \bb\label{t2} A=\frac{1}{2}\left [ -{\nabla }^s {\xi }^k -{\nabla }^k {\xi }^s
 + g^{ks} div (\xi )+2a g^{ks}
 \right ]\sqrt{g} u_k u_s .\ee
 From (\ref{t1}), (\ref{h44}) and (\ref{t2}) we obtain (\ref{n2}).

 \

 \rm Now we shall prove the first part of Theorem 2, namely

 \

 {\bf Lemma 6.} {\it For an arbitrary $f(u)$
 any symmetry of $(\ref{1.1})$ is a variational symmetry, that is,
 the isometry group of $M^{n}$ and the variational symmetry group
 of $(\ref{1.1})$ coincide.}

 \

 {\bf Proof.} \rm We have already seen in section 4 that in this
 case $a=b=\mu=div(\xi )=0$. Substituting this data into
 (\ref{n2}) we obtain $X^{(1)}L + L\;D_i{\xi }^i=0 $. Thus $X=\xi^{i}(x){\f{\p}{\p
 x^{i}}}$ is a variational symmetry. And vice-versa: from
 (\ref{n2}) it follows that any variational symmetry of (\ref{1.1})
 with arbitrary $f(u)$ is an isometry.

 \

 {\bf Lemma 7.} {\it In the exponential case $f(u)=e^u$ the only variational
 symmetries are the isometries of} $M$.

 \

 {\bf Proof.} \rm Substituting $a=0$, $b=-\mu =const$, $div(\xi )= n\mu /2$ and
 $L=g^{ks}u_ku_s/2-e^u$ into (\ref{n2}) we obtain
 \bb\label{n3} X^{(1)}L + L\;D_i{\xi }^i =  (n-2)\mu\; L/2.\ee
 Hence it is clear that $X$ is never a divergence symmetry. Again
 from (\ref{n3}), the symmetry $X$ is variational if and only if
 $n=2$ or $\mu =0$. Since $n\geq 3$ it follows that the only
 variational symmetries in the exponential case are the
 isometries ($\mu =0$).

 \

 {\bf Lemma 8.} {\it In the power case $f(u)=u^p$, $p\neq 0$, $p\neq 1$,
 the symmetry $(\ref{z10})$ is variational if and only if
 \[ p=\frac{n+2}{n-2},\]
 that is, $p+1$ equals to the critical Sobolev exponent.}

 \

 {\bf Proof.} \rm We put $a=\mu /(1-p)$, $\mu-$const, $b=0$, $div(\xi )= n\mu /2$ and
 $L=g^{ks}u_ku_s/2-u^{p+1}/(p+1)$ into (\ref{n2}). We obtain
 \bb\label{n4} X^{(1)}L + L\;D_i{\xi }^i
 =(-1+\f{n}{2}+\f{2}{1-p})\;\mu\;\sqrt{g}g^{ks}u_k u_s -
 (\f{n}{2}\f{1}{p+1}+\f{1}{1-p})\;\mu\;\sqrt{g}\;u^{p+1}. \ee
 Hence $X$ is a variational symmetry if and only if in (\ref{n4}) the
 coefficients of the terms containing $u$ and its derivatives
 vanish, which holds if and only if $p=(n+2)/(n-2)$.

 \

 {\bf Lemma 9.} {\it Let $X$ be the Lie point symmetry $(\ref{z12})$.
 Then}
 \bb\label{n5}X^{(1)}L + L\;D_i{\xi
 }^i=\f{2-n}{4}\sqrt{g}g^{ij}{\mu }_i\;u u_j.\ee

 \

 {\bf Proof.} \rm Substituting $a=(2-n)\mu /4$, $b=0$, $div(\xi )= n\mu /2$ and
 $L=g^{ks}u_ku_s/2-(n-2)u^{2n/(n-2)}/(2n)$ into (\ref{n2}) we
 obtain (\ref{n5}).

 \

 {\bf Lemma 10.} {\it If $\lb \mu =0$, the following equality holds:
 \bb\label{n6} \f{2-n}{4}\sqrt{g}g^{ij}{\mu }_j\;u u_i = D_i
 {\varphi }^i ,\ee
 where}
 \[ {\varphi }^i = {\varphi }^i (x,u) =\f{2-n}{8}\sqrt{g}g^{ij}{\mu
 }_j\;u^2 .\]

 \

 {\bf Proof.} \rm We calculate the total divergence of $\varphi = ({\varphi
 }^i)$:

 \[ \begin{array}{lll}D_i {\varphi }^i &= &(\f{\p}{\p x^{i}} + u_i \f{\p}{\p u} +
 u_{ik} \f{\p}{\p u_k}+...) (\f{2-n}{8}\sqrt{g}g^{ij}{\mu
 }_j\;u^2 ) \\
 & & \\
 &=& \f{2-n}{8}(\sqrt{g}g^{ij}{\mu }_j)_i\;u^2 + \f{2-n}{4}\sqrt{g}g^{ij}{\mu
 }_j\;u u_i  \\
 & & \\
 & = & \f{2-n}{8}[ {\nb }_i (\sqrt{g}g^{ij}{\mu }_j) -
 {\Ga }^i_{ki} \sqrt{g}g^{kj}{\mu }_j]\;u^2 + \f{2-n}{4}\sqrt{g}g^{ij}{\mu
 }_j\;u u_i  \\
 & & \\
 & = & \f{2-n}{8}[ {\nb }_i (\sqrt{g})g^{ij}{\mu }_j + \sqrt{g} {\nb }_i (g^{ij}{\mu }_j) -
 {\Ga }^i_{ki} \sqrt{g}g^{kj}{\mu }_j]\;u^2 + \f{2-n}{4}\sqrt{g}g^{ij}{\mu
 }_j\;u u_i  \\
 & & \\
 & = & \f{2-n}{8} [({\nb }_i (\sqrt{g})- {\Ga }^i_{ki}\sqrt{g})g^{ij}{\mu }_j +
 \sqrt{g}\lb \mu ]\;u^2
 + \f{2-n}{4}\sqrt{g}g^{ij}{\mu }_j\;u u_i  \\
 & & \\
 & = &  \f{2-n}{4}\sqrt{g}g^{ij}{\mu }_j\;u u_i \end{array} \]
 since the metric is parallel with respect to the Levi-Civita connection, $\lb \mu =0$ and
 ${\Ga }^i_{ki}\sqrt{g}= (\ln \sqrt{g})_i \sqrt{g}= (\sqrt{g})_i ={\nb
 }_i (\sqrt{g})$.

 \

 \rm Then from (\ref{n5}) and (\ref{n6}) it follows that \[X^{(1)}L + L\;D_i{\xi
 }^i= D_i {\varphi }^i ,\] that is, $X$ is a divergence symmetry.
 In this way we have proved the part 3.1) of Theorem 2.

 \

 {\bf Lemma 11.} {\it Let $X$ be the Lie point symmetry $(\ref{z15})$.
 Then}
 \bb\label{n7}X^{(1)}L + L\;D_i{\xi
 }^i=-\f{1}{2}\sqrt{g}\lb \mu\;u^2 - \sqrt{g}g^{ij}{\mu }_i\;u u_j +
 \f{1}{2} \sqrt{g}g^{ij}(\lb \mu )_i \;u_j.\ee

 \

 {\bf Proof.} \rm Substituting $a=- \mu $, $b=\lb\mu /2$, $div(\xi )= 3\mu $ and
 $L=g^{ks}u_ku_s/2-u^{3}/3$ into (\ref{n2}) we
 obtain (\ref{n7}).

 \

 {\bf Lemma 12.} {\it If ${\lb }^{\!\!\!2}\mu =0$, the following
 equality holds:
 \bb\label{n8}-\f{1}{2}\sqrt{g}\lb \mu\;u^2 - \sqrt{g}g^{ij}{\mu }_i\;u u_j +
 \f{1}{2} \sqrt{g}g^{ij}(\lb \mu )_i \;u_j= D_i
 {\phi }^i ,\ee
 where}
 \[ {\phi }^i = {\phi }^i (x,u) =-\f{1}{2}\sqrt{g}g^{ij}{\mu
 }_j\;u^2 + \sqrt{g}g^{ij}(\lb \mu )_j \;u.\]

 \

 {\bf Proof.} \rm We calculate the total divergence of $\phi = ({\phi
 }^i)$:

 \[ \begin{array}{lll}D_i {\phi }^i &= &(\f{\p}{\p x^{i}} + u_i \f{\p}{\p u} +
 u_{ik} \f{\p}{\p u_k}+...) (-\f{1}{2}\sqrt{g}g^{ij}{\mu
 }_j\;u^2 + \sqrt{g}g^{ij}(\lb \mu )_j \;u) \\
 & & \\
 &=& -\f{1}{2}(\sqrt{g}g^{ij}{\mu }_j)_i\;u^2 - \sqrt{g}g^{ij}{\mu
 }_j\;u u_i + (\sqrt{g}g^{ij}(\lb \mu )_j )_i\;u  + \sqrt{g}g^{ij}(\lb \mu )_j \;u_i\\
 & & \\
 & = & -\f{1}{2}[ {\nb }_i (\sqrt{g}g^{ij}{\mu }_j) -
 {\Ga }^k_{ik} \sqrt{g}g^{ij}{\mu }_j]\;u^2 - \sqrt{g}g^{ij}{\mu
 }_j\;u u_i \\
 & & \\
 & &+{\nb }_i(\sqrt{g}g^{ij}(\lb\mu
 )_j\;u  -  {\Ga }^k_{ik} \sqrt{g}g^{ij}(\lb \mu )_j\;u + \sqrt{g}g^{ij}(\lb \mu )_j \;u_i\\
 & & \\
 & = & -\f{1}{2}[({\nb }_i (\sqrt{g})- {\Ga }^k_{ik}\sqrt{g}]g^{ij}{\mu }_j\;u^2
 - \f{1}{2}\sqrt{g}\lb\mu\;u^2 - \sqrt{g}g^{ij}{\mu
 }_j\;u u_i \\
 & & \\
 & & +[({\nb }_i (\sqrt{g})- {\Ga }^k_{ik}\sqrt{g}]g^{ij} (\lb\mu )_j\;u + \sqrt{g} {\lb }^{\!\!\!2}\mu\;u + \sqrt{g}g^{ij} (\lb\mu )_i\;u_j \\
 & & \\
 & = & -\f{1}{2}\sqrt{g}\lb \mu\;u^2 -
 \sqrt{g}g^{ij}{\mu }_i\;u u_j + \f{1}{2} \sqrt{g}g^{ij}(\lb \mu )_i \;u_j \end{array} \]
 since the metric is parallel with respect to the Levi-Civita connection,
 ${\lb }^{\!\!\!2}\mu =0$ and
 ${\Ga }^i_{ki}\sqrt{g}= (\ln \sqrt{g})_i \sqrt{g}= (\sqrt{g})_i ={\nb
 }_i (\sqrt{g})$.

 \

 \rm Then from (\ref{n7}) and (\ref{n8}) it follows that \[X^{(1)}L + L\;D_i{\xi
 }^i= D_i {\phi }^i ,\] that is, $X$ is a divergence symmetry.
 In this way we have proved the part 3.2) of Theorem 2.

  \

 \section{The Noether symmetries in linear cases}

 \

 Here we shall prove the part 4 of Theorem 2.

\

{\bf Lemma 13.} {\it Let $X$ be the symmetry $(\ref{z2})$. Then}

\bb\label{n10}X^{(1)}L + L\;D_i{\xi
 }^i=2c\,L+\f{2-n}{4}\sqrt{g}g^{ij}{\mu }_i\;u u_j+\sqrt{g}g^{ij}b_{i}u_{j}.\ee
\rm

 {\bf Proof.} \rm Substituting $a=(\frac{2-n}{4}\;\mu (x)
 +c) $, $div(\xi )= \f{n}{2}\mu $, $F(u)=f(u)=0$ and
 $L=\sqrt{g} g^{ks}u_ku_s/2$ into (\ref{n2}) we
 obtain (\ref{n10}).

\

{\bf Lemma 14.} {\it If $\lb \mu=\lb b =0$, the following equality
holds: \bb\label{n11} \f{2-n}{4}\sqrt{g}g^{ij}{\mu }_j\;u u_i
+\sqrt{g}g^{ij}b_{i}u_{j}= D_i {\phi }^i ,\ee where}
\bb\label{n12} {\phi }^i = {\phi }^i (x,u)
=\f{2-n}{8}\sqrt{g}g^{ij}{\mu }_j\;u^2+\sqrt{g}g^{ij}b_{j}\;u .
\ee

{\bf Proof.} \rm We have $\phi^{i}=\varphi^{i}+\psi^{i}$, where
$\varphi^{i}=\f{2-n}{8}\sqrt{g}g^{ij}{\mu }_j\;u^2$ and $\psi^{i}
= \sqrt{g}g^{ij}b_{j}\;u$. Since the divergence of $\varphi^{i}$
was already calculated in Lemma 10 we have that

\[ \begin{array}{lll}\ds{D_i \phi^{i}}=D_i\varphi^{i}+D_i\psi^i &= &
\f{2-n}{4}\sqrt{g}g^{ij}{\mu }_i\;u u_j+\ds{(\f{\p}{\p x^{i}} +
u_i \f{\p}{\p u} +
u_{ik} \f{\p}{\p u_k}+...) (\sqrt{g}g^{ij}b_{j}\;u)} \\
& & \\
&=& \f{2-n}{4}\sqrt{g}g^{ij}{\mu }_i\;u
u_j+\ds{(\sqrt{g}g^{ij}b_j)_i\;u +
\sqrt{g}g^{ij}b_j\; u_i}  \\
& & \\
& = & \f{2-n}{4}\sqrt{g}g^{ij}{\mu }_i\;u u_j+\ds{\sqrt{g}\lb
b\;u+\sqrt{g}g^{ij}b_{i}u_{j}} \\ & & \\ &=&
\f{2-n}{4}\sqrt{g}g^{ij}{\mu }_i\;u u_j+\sqrt{g}g^{ij}b_{i}u_{j}
\end{array} \]
since $\lb b =0$.

\

{\bf Lemma 15.}  {\it Let $X$ be the symmetry $(\ref{z2})$. Then
$X$ is a Noether symmetry if and only if} $c=0$.

\

{\bf Proof.} \rm  From (\ref{n10}) and (\ref{n11}) it follows that
\bb\label{n13}X^{(1)}L + L\;D_i{\xi }^i=2c\;L +D_i {\phi }^i. \ee
Then the conclusion of Lemma 15 follows from (\ref{n13}). This
proves part 4.1) of Theorem 2.

\

We observe that the symmetry $(\ref{z2})$ can be written as
$$X=\xi+\frac{2-n}{4}\;\mu (x)\f{\p }{\p u}
 +b(x)\frac{\p}{\p u}+c\; u\f{\p }{\p u} ,
$$
where $\xi=\xi^{i}(x)\f{\p}{\p x^{i}}$ is a conformal Killing
vector field satisfying (\ref{z5}), $c=const.$, $\lb b=0$ and
$\lb\mu=0$. The potentials $\varphi^{i},\,\psi^{i}$ are the
potentials of the symmetries $\xi+\frac{2-n}{4}\;\mu (x)\f{\p }{\p
u}$ and $b(x)\frac{\p}{\p u},$ respectively. The symmetry $u\f{\p
}{\p u}$ corresponds to a non-Noetherian symmetry and it reflects
the linearity of the equation.

\

{\bf Lemma 16.} {\it Let $X$ be the symmetry $(\ref{h1})$. Then}
\bb\label{n14}X^{(1)}L + L\;D_i{\xi
 }^i=2c\,L+\f{2-n}{4}\sqrt{g}g^{ij}{\mu }_i\;u u_j+
 \sqrt{g}g^{ij}b_{i}u_{j}+\sqrt{g}\lb b \;u-\sqrt{g}\;b\;k.\ee
\rm

{\bf Proof.} \rm Substituting \[a=(\frac{2-n}{4}\;\mu
(x)+c),\,\,div(\xi )= \f{n}{2}\mu,\,\,F(u)=ku,\,f(u)=k,\,
L=\sqrt{g} g^{ks}u_ku_s/2-ku\sqrt{g} \] into (\ref{n2}) and using
$\lb b=ck-\f{n+2}{4}\mu k$ (see (\ref{h2})) we obtain (\ref{n14}).

\

Observe that by (\ref{h2}): $\mu =\f{4}{(n+2)}(c-\f{1}{k}\lb b)$.
Hence $\lb\mu=0$ since $ {\lb }^{\!\!\!2} b =0 $.

\

{\bf Lemma 17.} {\it Let $\phi^{i}$ be the potential $(\ref{n12})$
and $\lb\mu=0$. Then} \bb\label{n15}
D_{i}\phi^{i}=\f{2-n}{4}\sqrt{g}g^{ij}{\mu }_i\;u
u_j+\sqrt{g}g^{ij}b_{i}u_{j}+\sqrt{g} \lb b\;u \ee

\

\rm The proof of this lemma is similar to that of Lemma 10 and
Lemma 14.

\

{\bf Lemma 18.} {\it Let $X$ be the symmetry $(\ref{h1})$. Then
$X$ is a Noether symmetry if and only if} $c=b=0$.

\

{\bf Proof.} \rm Substituting $\lb b=-\f{n+2}{4}k\mu+kc$ and
(\ref{n15}) into (\ref{n14}), we obtain \bb\label{n16}X^{(1)}L +
L\;D_i{\xi
 }^i=2c\,L+D_{i}\phi^{i}-\sqrt{g}\;b\;k \ee
 which implies Lemma 18. This proves part 4.2) of Theorem 2.

\

Now let $X$ be the symmetry (\ref{z2}) of (\ref{1.1}) with
$f(u)=u$. We shall only sketch the proof of 4.3) of Theorem 2.

Substituting $a=(2-n)\mu /4 +c$, $div(\xi )=n/2$, $f(u)=u$,
$F(u)=u^2/2$, $L=\sqrt{g}g^{ij}u_iu_j/2 - \sqrt{g} u^2/2$ into
(\ref{n2}), we obtain \bb\label{wv1}X^{(1)}L + L\;D_i{\xi
 }^i=2c\,L+D_{i}\phi^{i} -\sqrt{g} u^2 (\frac{2-n}{2}\lb \mu
 +\mu)/2 - \sqrt{g} u (\lb b +b), \ee
 where $\phi^i $ is given in (\ref{n12}). Then from (\ref{z6}),
 (\ref{z7})and (\ref{wv1}) we get
 \[ X^{(1)}L + L\;D_i{\xi
 }^i=2c\,L+D_{i}\phi^{i}.\]
 Hence $X$ is a Noether symmetry if and only if $c=0$.

Thus, we have concluded the proof of Theorem 2.

\

\section{Examples: Poisson equations on Thurston geometries}

\

In this section we apply our results to Poisson equations on the
Thurston geometries. The presentation is very schematic in order
not to increase the volume of the paper.

All examples presented here involve elliptic forms of equation
(\ref{1.1}). Examples of symmetry analysis involving some
particular hyperbolic cases of (\ref{1.1}) can be found in
\cite{A, igor,fu1}.

Some of the results presented in this section were verified using
the SYM package \cite{grego1, grego2}.

To the authors' knowledge, the results in \ref{hyperbolic},
\ref{sphere}, \ref{sol}, \ref{prod1}, \ref{prod2}, \ref{sl2} and
\ref{heisenberg} are original.

\

\subsection{Thurston geometries}

\

A manifold $M^{n}$ is said to be homogeneous if, for every
$x,\;y\in M$, there exists an isometry of $M^{n}$ such that it
leaves $x$ in $y$. Let $X$ be the universal covering of $M^{n}$
and $G$ its isometry group.

A geometry consists of a pair $(X,G)$ as above, where $X$ is a
connected manifold and $G$ is a group that acts effectively and
transitively on $X$, and where all stabilizers $G_{x}$ are
compact. This is also equivalent to the data of a connected Lie
group $G$ and a compact Lie subgroup $H$ of $G$, if we associate
to this data the homogeneous space $X = G/H$ endowed with the
natural left action of $G$.

Two geometries $(X,G)$ and $(X',G')$ are identified if there is a
diffeomorphism from $X$ to $X'$ which sends the action of $G$ to
the action of $G'$. $(X,G)$ is said to be maximal if there is no
larger geometry $(X',G')$ with $G \subseteq G'$ and $G \neq G'$.
For more details, see \cite{bo,scott}.

There are exactly 8 three-dimensional maximal geometries $(X,G)$,
the so-called Thurston geometries.

Thurston \cite{Tu1} has classified the three-dimensional,
simply-connected, homogeneous manifolds as follows (see also
\cite{Tu, Tu2, bo, scott}):

\begin{itemize}
\item the Euclidean space $\R^3=\{(x,y,z)\,|\,x,y,z\in\R\}$, with
canonical metric $$ds^{2}=dx^2+dy^2+dz^2;$$ \item the Hyperbolic
space $\mathbb{H}^3=\{(x,y,z)\in\R^3\,|\,z>0\}$, with metric
$$ds^{2}=(dx^2+dy^2+dz^2)/z^2;$$ \item the Sphere
$\mathbb{S}^3=\{(x_1,x_2,x_3,x_4)\in\R^4\,|\,\sum_{i=1}^4
x_i^2=1\}$, with induced metric from $\R^{4}$; \item the solvable
group Sol, which can be defined as the Lie group $(\R^3,\,\ast)$
where $$(x,y,z)*(x',y',z')=(x+x',y+e^{-x}\,y',z+e^x\,z'),$$  with
left-invariant metric $ds^{2}=dx^2+e^{2x}\,dy^2+e^{-2x}\,dz^2$;
\item the space $\mathbb{S}^2\times\R$, with product metric; \item
the space $\mathbb{H}^2\times\R$, with product metric. Here
$\mathbb{H}^2$ is the two-dimensional Hyperbolic space; \item the
universal covering of ${SL_2}(\R)$, or
$\R^3_+=\{(x,y,z)\in\R^3\,|\,z>0\}$, with metric
\[
\ds{ds^{2}=\bigg(dx+\frac{dy}{z}\bigg)^2+\frac{(dy^2+dz^2)}{z^2}}\]
and \item the Heisenberg group $H^{1}$, whose group structure is
given by \[\phi((x,y,t),  (x',y',t')) = (x+x', y+y', t+t' +2(y
x'-xy'))\] and the left-invariant metric is
$ds^{2}=dx^2+dy^2+\big(dz+2ydx-2xdy\big)^2$.
\end{itemize}

Three of them are isotropic geometries: if the curvature is
positive, then the isotropic geometry is the 3-Sphere
$\mathbb{S}^3$. If the curvature is negative, the isotropic
geometry is the Hyperbolic space $ \mathbb{H}^3$. If the curvature
vanishes, the isotropic geometry is the Euclidean space $\R^{3}$.

Four of the Thurston geometries, the product spaces
$\mathbb{S}^2\times\R$, $\mathbb{H}^2\times\R$ and two Lie groups,
$\widetilde{SL_2}(\R)$ and $H^{1}$, are known as the four Seifert
type geometries.

Finally, we have the Sol group, which possesses this name because
the group $G$ of the pair $(Sol, G)$ is solvable and it is the
only one of the Thurston geometries with this property.

For more details about the Thurston geometries, see \cite{Tu, Tu1,
Tu2, bo, scott}.

\

\subsection{The Euclidean Space}\label{euclidean}

 \

The application presented in this section is well-known. It
corresponds to the group classification, Noether symmetries and
conservation laws of nonlinear Poisson equations in $\R^{3}$. The
group classification of these equations can be found in \cite{sv}
as a particular case. The Noether symmetries and conservation laws
are established in \cite{yb2} in a more general context.

Here we shall consider the three-dimensional vector space $\R^{3}$
with the Euclidean metric $$ds^{2}=dx^{2}+dy^{2}+dz^{2}$$ and the
Poisson equation \bb\label{r1} u_{xx}+u_{yy}+u_{zz}+f(u)=0. \ee

\subsubsection{The group classification}

\begin{enumerate}
\item For any arbitrary function $f(u)$, the symmetry group of
(\ref{r1}) coincides with the isometry group of $\R^{3}$. It is
well-known (see e.g. \cite{i,dnf}) that the latter is generated by
translations and rotations given by \bb\label{e1} \ba{l}
\ds{R_{1}=\f{\p}{\p x}},\,\,\,\ds{R_{2}=\f{\p}{\p y}},\,\,\,\ds{R_{3}=\f{\p}{\p z}},\\
\\
\ds{R_{4}=y\f{\p}{\p x}-x\f{\p}{\p y}},\,\,\,\ds{R_{5}=y\f{\p}{\p
z}-z\f{\p}{\p y}}, \,\,\,\ds{R_{6}=z\f{\p}{\p x}-x\f{\p}{\p z}}.
\ea \ee Hence (\ref{e1}) determine the symmetry group of
(\ref{r1}) for arbitrary function $f(u)$.

\item If $f(u)=0$ then the additional to (\ref{e1}) symmetries are
\bb\label{e2} \ba{l}
\ds{R_{7}=x\f{\p}{\p x}+y\f{\p}{\p y}+z\f{\p}{\p z} +\f{u}{2} \frac{\p}{\p u}},\\
\\
\ds{R_{8}=xz\f{\p}{\p x}+yz\f{\p}{\p y}+\f{(z^{2}-x^{2}-y^{2})}{2}\f{\p}{\p z} -zu \frac{\p}{\p u}},\\
\\
\ds{R_{9}=xy\f{\p}{\p x}+\f{y^{2}-x^{2}-z^{2}}{2}\f{\p}{\p y}+yz\f{\p}{\p z} -\f{yu}{2}\frac{\p}{\p u}},\\
\\
\ds{R_{10}=\f{x^{2}-y^{2}-z^{2}}{2}\f{\p}{\p x}+xy\f{\p}{\p y}+xz\f{\p}{\p z} -\f{xu}{2} \frac{\p}{\p u}}.\\
\\
\ea \ee \bb\label{e3} R_{11}=u\f{\p}{\p
u},\,\,\,\,\,\,R_{\infty}=b(x,y,z)\f{\p}{\p u}, \ee where $\Delta
b=0$.

\item The case $f(u)=k=const\neq 0$ reduces to the homogeneous
case under the change $u\rightarrow u-k\,x^{2}/2$.

\item If the function $f$ is a linear function, $f(u)=u$, then the
additional symmetry generators are given by $(\ref{e3})$, with
$\Delta b+b=0$

\item For exponential nonlinearity $f(u)=e^u$, the additional
generator is
$$R_{13}=x\f{\p}{\p x}+y\f{\p}{\p y}+z\f{\p}{\p z} -2 \frac{\p}{\p u}.$$

\item For power nonlinearity $f(u)=u^p$, $p\neq 0$, $p\neq 1$, and
$p\neq 5$, the additional generator is
$$R_{14}=x\f{\p}{\p x}+y\f{\p}{\p y}+z\f{\p}{\p z} +\f{2}{1-p}u \frac{\p}{\p u}.$$

\item If $p=5$, then the additional infinitesimal generators of
the Lie point symmetries are given in $(\ref{e2})$.

\end{enumerate}

We observe that the critical Sobolev exponent $\f{n+2}{n-2}$ in
this case is exactly 5 and, if $p(p-1)(p-5)\neq 0$, then the
symmetry group of (\ref{r1}), with $f(u)=u^{p}$, is isomorphic to
the symmetry group of (\ref{r1}) with $f(u)=e^{u}$ and the special
conformal group generated by the symmetries is the group of
homothetic motions in $\R^{3}$.

\subsubsection{The Noether symmetries}

\begin{enumerate}
\item The isometry group is a variational symmetry group of the
nonlinear Poisson equation in $\R^{3}$. In particular, it is the
Noether symmetry group of the cases $f(u)=e^{u}$ and $f(u)=u^{p}$,
with $p\neq 0,\,1,\,5$.

\item The conformal group of $\R^{3}$ and the symmetry
$R_{\infty}$ generate a Noether symmetry group for $\Delta u=0$.

\item The isometry group and the symmetry $R_{\infty}$ generate a
Noether symmetry group for $\Delta u+u=0$.

\item The full conformal group of $\R^{3}$ is a Noether symmetry
group for $\Delta u+u^{5}=0$.

\end{enumerate}

\rm

The corresponding conservation laws can be obtained as a
particular case of those, more general, established in \cite{yb2}
and for this reason they are not stated explicitly here.

\

\subsection{The Hyperbolic Space}\label{hyperbolic}

\

We consider the Klein's model of the Hyperbolic space
$\mathbb{H}^3$ represented by the set of $(x,y,z)\in\R^{3}$, with
$z>0$, and endowed with the metric
$$ds^{2}=\f{dx^{2}+dy^{2}+dz^{2}}{z^{2}}.$$

This metric has constant negative scalar curvature $R=-6$ and its
sectional curvature is equal to $-1$ (see \cite{carmo}, p. 160).
Thus one immediately concludes from \cite{ya}, p. 57, that the
isometry group of $\mathbb{H}^3$ possesses a 6 dimensional Lie
algebra.

It is easy to check that the following vector fields

\bb\label{hy1} \ba{lcl}
\ds{H_{1}=\f{\p}{\p x}}, \,\, \ds{H_{2}=\f{\p}{\p y}}, \,\, \ds{H_{3}=-y\f{\p}{\p x}+x\f{\p}{\p y}}, & & \\
\\
\ds{H_{4}=x\f{\p}{\p x}+y\f{\p}{\p y}+z\f{\p}{\p z}}, & &\\
\\
\ds{H_{5}=\f{x^{2}-y^{2}-z^{2}}{2}\f{\p}{\p x}+xy\f{\p}{\p y}+xz\f{\p}{\p z}}, & &\\
\\
\ds{H_{6}=xy\f{\p}{\p x}+\f{-x^{2}+y^{2}-z^{2}}{2}\f{\p}{\p
y}}+yz\f{\p}{\p z}& & \ea \ee are Killing fields on
$(\mathbb{H}^3,g)$ and, for maximality, they form a basis of
generators of $Isom(\mathbb{H}^3,g)$.

The nonlinear Poisson equation on $(\mathbb{H}^3,g)$ is given by
\bb\label{hy2} z^{2}\Delta u-zu_{z}+f(u)=0, \ee where $\Delta$ is
the Laplace operator on $\R^{3}$ and $u_z=\f{\kd u}{\kd z}$.

\subsubsection{The group classification}

\begin{enumerate}
\item For any function $f(u)$, the symmetry group coincides with
$Isom(\mathbb{H}^3,g)$.

\item  If $f(u)=0$ then the additional symmetries are
\bb\label{hy3} H_{\infty}=b\f{\p}{\p u}, \ee with $\lb b=0$, and
\bb\label{hy4} H_{7}=u\f{\p}{\p u}. \ee

\item The case $f(u)=k=const.\neq 0$ is equivalent to case
$f(u)=0$ under the change $v=u +(k/2)\ln\,z$.

\item If the function $f$ is a linear function, $f(u)=u$, then the
additional symmetry generator is given by $(\ref{hy3})$, with
$\Delta b+b=0$, and (\ref{hy4}).
\end{enumerate}

\subsubsection{The Noether symmetries}

\begin{enumerate}

\item The isometry group of $(\mathbb{H}^3,g)$ is a variational
symmetry group. \item Symmetries $(\ref{hy3})$, with $\lb b=0$ or
$\lb b+b=0$, are the Noether symmetries to the cases $f(u)=0$ or
$f(u)=u$, respectively.
\end{enumerate}

\subsubsection{The conservation laws}

Here we present the conservation laws corresponding to the Noether
symmetries of the equations (\ref{hy2}) with arbitrary $f(u)$.

\

\begin{enumerate}
\rm\item For the symmetry $H_{1}$, the conservation law is
$Div(A)=0$, where $A=(A_{1},A_{2},A_{3})$ and
$$
\ba{l c l}
A_{1} & = & \ds{\f{y^{2}+z^{2}-x^{2}}{4z}\left(u_{x}^{2}-u_{y}^{2}-u_{z}^{2}\right)-\f{1}{z}\left(xyu_{x}u_{y}+xzu_{x}u_{z}\right)+\f{y^{2}+z^{2}-x^{2}}{4z^{2}}F(u)},\\
\\
A_{2} & = & \ds{\f{xy}{2z}\left(u_{x}^{2}-u_{y}^{2}+u_{z}^{2}\right)+\f{y^{2}+z^{2}-x^{2}}{2z}u_{x}u_{y} -xzu_{y}u_{z}-\f{xy}{2z^{2}}F(u)},\\
\\
A_{3} & = & \ds{\f{x}{2}\left(u_{x}^{2}+u_{y}^{2}-u_{z}^{2}\right)\f{y^{2}+z^{2}-x^{2}}{2z}u_{x}u_{z} -xyu_{y}u_{z}-\f{x}{2z}F(u)}.\\
\ea
$$

\rm\item For the symmetry $H_{2}$, the conservation law is
$Div(B)=0$, where $B=(B_{1},B_{2},B_{3})$ and
$$
\ba{l c l}
B_{1} & = & \ds{\f{xy}{2z}\left(u_{y}^{2}+u_{z}^{2}-u_{x}^{2}\right)+\f{x^{2}-y^{2}+z^{2}}{2z}u_{x}u_{y}-yu_{x}u_{z}-\f{xy}{2z^{2}}F(u)},\\
\\
B_{2} & = & \ds{\f{x^{2}+y^{2}+z^{2}}{2z}\left(u_{y}^{2}-u_{x}^{2}-u_{z}^{2}\right)-\f{xy}{z}u_{x}u_{y}-yu_{y}u_{z}-\f{1}{4z^{2}}F(u)},\\
\\
B_{3} & = & \ds{\f{y}{2}\left(u_{x}^{2}+u_{y}^{2}-u_{z}^{2}\right)-\f{xy}{z}u_{x}u_{y}+\f{x^{2}-y^{2}+z^{2}}{2z}u_{y}u_{z}-\f{y}{2z}F(u)}.\\
\ea
$$

\rm\item For the symmetry $H_{3}$, the conservation law is
$Div(C)=0$, where $C=(C_{1},C_{2},C_{3})$ and
$$
\ba{l c l}
C_{1} & = & \ds{\f{x}{2z}\left(u_{y}^{2}+u_{z}^{2}-u_{x}^{2}\right)-\f{y}{z}u_{x}u_{y}-u_{x}u_{z}-\f{x}{2z^{2}}F(u)},\\
\\
C_{2} & = & \ds{\f{y}{2z}\left(u_{x}^{2}-u_{y}^{2}+u_{z}^{2}\right)-\f{x}{z}u_{x}u_{y}-u_{y}u_{z}-\f{y}{2z^{2}}F(u)},\\
\\
C_{3} & = & \ds{\f{1}{2}\left(u_{x}^{2}+u_{y}^{2}-u_{z}^{2}\right)-\f{x}{z}u_{x}u_{z}-\f{y}{z}u_{y}u_{z}-\f{1}{2z}F(u)}.\\
\ea
$$

\rm\item For the symmetry $H_{4}$, the conservation law is
$Div(D)=0$, where $D=(D_{1},D_{2},D_{3})$ and
$$
\ba{l c l}
D_{1} & = & \ds{\f{y}{2z}\left(u_{y}^{2}+u_{z}^{2}-u_{x}^{2}\right)+\f{x}{z}u_{x}u_{y}-\f{y}{2z^{2}}F(u)},\\
\\
D_{2} & = & \ds{\f{x}{2z}\left(u_{y}^{2}-u_{x}^{2}-u_{z}^{2}\right)-\f{y}{z}u_{x}u_{y}+\f{1}{2z^{2}}F(u)},\\
\\
D_{3} & = & \ds{\f{x}{z}u_{y}u_{z}-\f{y}{z}u_{x}u_{z}}.\\
\ea
$$

\rm\item For the symmetry $H_{5}$, the conservation law is
$Div(E)=0$, where $E=(E_{1},E_{2},E_{3})$ and
$$
\ba{l c l}
E_{1} & = & \ds{\f{1}{2z}\left(u_{y}^{2}+u_{z}^{2}-u_{x}^{2}\right)-\f{1}{2z^{2}}F(u)},\\
\\
E_{2} & = & \ds{-\f{1}{z}u_{x}u_{y}},\\
\\
E_{3} & = & \ds{-\f{1}{z}u_{x}u_{z}}.\\
\ea
$$

\rm\item For the symmetry $H_{6}$, the conservation law is
$Div(F)=0$, where $F=(F_{1},F_{2},F_{3})$ and
$$
\ba{l c l}
F_{1} & = & \ds{-\f{1}{z}u_{x}u_{y}},\\
\\
F_{2} & = & \ds{\f{1}{2z}\left(u_{x}^{2}-u_{y}^{2}+u_{z}^{2}\right)-\f{1}{2z^{2}}F(u)},\\
\\
F_{3} & = & \ds{-\f{1}{z}u_{y}u_{z}}.\\
\ea
$$

\rm\item For the symmetry $H_{\infty}$, with $\lb b=0$, the
conservation law is $Div(G)=0$, where $G=(G_{1},G_{2},G_{3})$ and
\bb\label{hy29}\ba{lcl}
G_{1} &= &\ds{\f{bu_{x}-b_{x}u}{z}},\\
\\
G_{2} & =&  \ds{\f{bu_{y}-b_{y}u}{z}},\\
\\
G_{3} & = & \ds{\f{bu_{z}-b_{z}u}{z}}. \ea \ee

\rm\item For the symmetry $H_{\infty}$, with $\lb b+b=0$, the
conservation law is $Div(G)=0$, where $G$ is given in
$(\ref{hy29})$.
\end{enumerate}

\

\subsection{The sphere}\label{sphere}

\

\rm Let us now consider the 3-sphere $\mathbb{S}^3$. Its metric is
given by the restriction to $\mathbb{S}^{3}$ of the canonical
metric of $\R^{4}$. Or, more specifically \bb\label{s1}
ds^{2}=\f{4}{(1+x^{2}+y^{2}+z^{2})^{2}}(dx^{2}+dy^{2}+dz^{2}). \ee

This metric determines the following Poisson equation on
$\mathbb{S}^3$ \bb\label{s4} \Delta
u-\f{2}{1+x^{2}+y^{2}+z^{2}}(xu_{x}+yu_{y}+zu_{z})+f(u)=0, \ee
where $\Delta$ denotes the Laplacian in $\R^{3}$ 

The $Isom(\mathbb{S}^{3},g)$ is generated by the following vector
fields: \bb\label{s5} \ba{l c l}
S_{1} & =& \ds{(1+x^{2}-y^{2}-z^{2})\f{\p}{\p x}+2xy\f{\p}{\p y}+2xz\f{\p}{\p z}},\\
\\
S_{2} &= & \ds{2xy\f{\p}{\p x}+(1-x^{2}+y^{2}-z^{2})\f{\p}{\p y}}+2zy\f{\p}{\p z},\\
\\
S_{3} &= & \ds{2xz\f{\p}{\p x}+2yz\f{\p}{\p y}}+(1-x^{2}-y^{2}+z^{2})\f{\p}{\p z},\\
\\
S_{4} &= & y\; \ds{\frac{\kd }{\kd x}-x\; \frac{\kd }{\kd y}}
\\
\\
S_{5} &= & z\; \ds{\frac{\kd }{\kd x}-x\; \frac{\kd }{\kd z}}\\
\\
S_{6} &= & z\; \ds{\frac{\kd }{\kd y}-y\; \frac{\kd }{\kd z}}. \ea
\ee

The scalar curvature of $(\mathbb{S}^3,g)$ is $R=6$.

\subsubsection{Group classification}

\begin{enumerate}
\item {\bf Arbitrary $f(u)$:} It is immediate that the vector
fields (\ref{s5}) are symmetries. (See Theorem 1 and Corollary 1.)

\item {\bf Linear case:} In addition to $Isom(\mathbb{S}^{3},g)$,
we have the symmetries \bb\label{s6} S_{7}=u\f{\p}{\p u} \ee and
\bb\label{s7} S_{\infty}=b\f{\p}{\p u}, \ee where $\lb b+b=0$,
$\int_{M^n} b\,dV=0$.

\item {\bf Homogeneous case:} In this case, the symmetries are
 given by (\ref{s5}), (\ref{s6}) and $$S_{8}=\frac{\p}{\p u}.$$
\end{enumerate}

\subsubsection{The Noether symmetries}

\begin{enumerate}
\item For arbitrary $f(u)$, the isometry group of
$(\mathbb{S}^{3},g)$ is a variational symmetry group. \item If
$f(u)=0$, in addition to the variational symmetries
$Isom(\mathbb{S}^{3},g)$, we have the divergence symmetry
$\frac{\p}{\p u}$. \item If $f(u)=u$, the additional divergence
symmetry is (\ref{s7}), with $\lb b+b=0$.
\end{enumerate}

\subsubsection{The conservation laws}

\begin{enumerate}
\item For the symmetry $S_{1}$, with arbitrary $f(u)$, the
conservation law is $Div(A)=0$, where $A=(A_{1},A_{2},A_{3})$ and
$$
\ba{l c l}
A_{1} & = & \ds{\f{(1+x^{2}-y^{2}-z^{2})(u_{y}^{2}+u_{z}^{2}-u_{x}^{2})-4xyu_{x}u_{y}+4xzu_{x}u_{z}}{(1+x^{2}+y^{2}+z^{2})^{2}}}\\
\\
& & \ds{-4\f{1+x^{2}-y^{2}-z^{2}}{(1+x^{2}+y^{2}+z^{2})^{3}}F(u)},\\
\\
A_{2} & = & \ds{\f{2(xyu_{x}^{2}-xyu_{y}^{2}+xzu_{z}^{2})-(1+x^{2}-y^{2}-z^{2})u_{x}u_{y}}{(1+x^{2}+y^{2}+z^{2})^{2}}}\\
\\
&&\ds{-\f{4xy}{(1+x^{2}+y^{2}+z^{2})^{3}}F(u)},\\
\\
A_{3} & = & \ds{\f{2((xzu_{x}^{2}+xzu_{y}^{2}-xzu_{z}^{2})-(1+x^{2}-y^{2}-z^{2})u_{x}u_{y})}{(1+x^{2}+y^{2}+z^{2})^{2}}}\\
\\
&&\ds{-\f{8xy}{(1+x^{2}+y^{2}+z^{2})^{3}}F(u)}.\\
\ea
$$

\item For the symmetry $S_{2}$, with arbitrary $f(u)$, the
conservation law is $Div(B)=0$, where $B=(B_{1},B_{2},B_{3})$ and
$$
\ba{l c l}
B_{1} & = & \ds{\f{2(xy(u_{y}^{2}+u_{z}^{2}-u_{x}^{2})-2yzu_{x}u_{z}-(1-x^{2}-y^{2}+z^{2})u_{x}u_{y})}{(1+x^{2}+y^{2}+z^{2})^{2}}}\\
\\
& &\ds{-\f{8xy}{(1+x^{2}+y^{2}+z^{2})^{3}}F(u)},\\
\\
B_{2} & = & \ds{\f{(1-x^{2}+y^{2}-z^{2})(u_{x}^{2}-u_{y}^{2}+u_{z}^{2})-4xyu_{x}u_{y}-4yzu_{y}u_{z}}{(1+x^{2}+y^{2}+z^{2})^{2}}}\\
\\
&
&+\ds{\f{4(1-x^{2}+y^{2}-z^{2})}{(1+x^{2}+y^{2}+z^{2})^{3}}F(u)},
\ea
$$
$$
\ba{lcl}
B_{3} & = & \ds{\f{2(yz(u_{x}^{2}+u_{y}^{2}-u_{z}^{2})-2xyu_{x}u_{z}-(1-x^{2}+y^{2}-z^{2})u_{x}u_{z}}{(1+x^{2}+y^{2}+z^{2})^{2}}}\\
\\
& &-\ds{\f{8yz}{(1+x^{2}+y^{2}+z^{2})^{3}}F(u)}.\\
\ea
$$

\item For the symmetry $S_{3}$, with arbitrary $f(u)$, the
conservation law is $Div(C)=0$, where $C=(C_{1},C_{2},C_{3})$ and
$$
\ba{l c l}
C_{1} & = & \ds{\f{2(xz(u_{y}^{2}+u_{z}^{2}-u_{x}^{2})-2yzu_{x}u_{y}-(1-x^{2}-y^{2}+z^{2})u_{x}u_{z})}{(1+x^{2}+y^{2}+z^{2})^{2}}}\\
\\
& &\ds{-\f{8xz}{(1+x^{2}+y^{2}+z^{2})^{3}}F(u)},\\
\\
C_{2} & = & \ds{\f{2(yz(u_{x}^{2}-u_{y}^{2}+u_{z}^{2})-2xzu_{x}u_{y}-(1-x^{2}-y^{2}+z^{2})u_{y}u_{z})}{(1+x^{2}+y^{2}+z^{2})^{2}}}\\
\\
& &\ds{-\f{8yz}{(1+x^{2}+y^{2}+z^{2})^{3}}F(u)}\\
\\
C_{3} & = & \ds{\f{(1-x^{2}-y^{2}+z^{2})(u_{x}^{2}+u_{y}^{2}-u_{z}^{2})-4xzu_{x}u_{z}-4yzu_{y}u_{z})}{(1+x^{2}+y^{2}+z^{2})^{2}}}\\
\\
& &\ds{-\f{4(1-x^{2}-y^{2}+z^{2})}{(1+x^{2}+y^{2}+z^{2})^{3}}F(u)}.\\
\ea
$$

\item For the symmetry $S_{4}$, with arbitrary $f(u)$, the
conservation law is $Div(D)=0$, where $D=(D_{1},D_{2},D_{3})$ and
$$
\ba{l c l}
D_{1} & = & \ds{\f{2y(u_{y}^{2}+u_{z}^{2}-u_{x}^{2})-4xu_{x}u_{y}}{(1+x^{2}+y^{2}+z^{2})^{2}}-\f{8y}{(1+x^{2}+y^{2}+z^{2})^{3}}F(u)},\\
\\
D_{2} & = & \ds{\f{-2x(u_{x}^{2}-u_{y}^{2}+u_{z}^{2})-4yu_{x}u_{y}}{(1+x^{2}+y^{2}+z^{2})^{2}}+\f{8x}{(1+x^{2}+y^{2}+z^{2})^{3}}F(u)}\\
\\
D_{3} & = & \ds{\f{4xu_{y}u_{z}-4yu_{x}u_{z}}{(1+x^{2}+y^{2}+z^{2})^{2}}}.\\
\ea
$$

\item For the symmetry $S_{5}$, with arbitrary $f(u)$, the
conservation law is $Div(E)=0$, where $E=(E_{1},E_{2},E_{3})$ and
$$
\ba{l c l}
E_{1} & = & \ds{\f{2z(u_{y}^{2}+u_{z}^{2}-u_{x}^{2})+4xu_{x}u_{z}}{(1+x^{2}+y^{2}+z^{2})^{2}}-\f{8z}{(1+x^{2}+y^{2}+z^{2})^{3}}F(u)},\\
\\
E_{2} & = & \ds{\f{4xu_{y}u_{z}-4zu_{x}u_{y}}{(1+x^{2}+y^{2}+z^{2})^{2}}}\\
\\
E_{3} & = & \ds{\f{2x(u_{x}^{2}+u_{y}^{2}-u_{z}^{2})+4zu_{x}u_{z}-4yu_{x}u_{z}}{(1+x^{2}+y^{2}+z^{2})^{2}}}-\f{8x}{(1+x^{2}+y^{2}+z^{2})^{3}}F(u).\\
\ea
$$

\item For the symmetry $S_{6}$, with arbitrary $f(u)$, the
conservation law is $Div(F)=0$, where $F=(F_{1},F_{2},F_{3})$ and
$$
\ba{l c l}
F_{1} & = & \ds{\f{4yu_{x}u_{z}-4zu_{x}u_{y}}{(1+x^{2}+y^{2}+z^{2})^{2}}},\\
\\
F_{2} & = & \ds{\f{2z(u_{x}^{2}+u_{y}^{2}-u_{z}^{2})+4yu_{y}u_{z}}{(1+x^{2}+y^{2}+z^{2})^{2}}-\f{8z}{(1+x^{2}+y^{2}+z^{2})^{3}}F(u)}\\
\\
F_{3} & = & \ds{\f{2y(u_{z}^{2}-u_{x}^{2}-u_{y}^{2})-4zu_{y}u_{z}}{(1+x^{2}+y^{2}+z^{2})^{2}}-\f{8y}{(1+x^{2}+y^{2}+z^{2})^{3}}F(u)}.\\
\ea
$$

\item For the symmetry $S_{\infty}$, with $\lb b+b=0$ or $\lb
b=0$, the conservation law is $Div(G)=0$, where
$G=(G_{1},G_{2},G_{3})$ and
$$ \ba{l c l}
G_{1} & = & \ds{\f{bu_{x}-b_{x}u}{1+x^{2}+y^{2}+z^{2}}},\\
\\
G_{2} & = & \ds{\f{bu_{y}-b_{y}u}{1+x^{2}+y^{2}+z^{2}}},\\
\\
G_{3} & = & \ds{\f{bu_{z}-b_{z}u}{1+x^{2}+y^{2}+z^{2}}}.\\
\ea
$$

\item For the symmetry $S_{8}=\f{\p}{\p u}$, the conservation law
is $Div(J)=0$, where $J=(J_{1},J_{2},J_{3})$ and
$$ \ba{l c l}
G_{1} & = & \ds{\f{u_{x}}{1+x^{2}+y^{2}+z^{2}}},\\
\\
G_{2} & = & \ds{\f{u_{y}}{1+x^{2}+y^{2}+z^{2}}},\\
\\
G_{3} & = & \ds{\f{u_{z}}{1+x^{2}+y^{2}+z^{2}}}.\\
\ea
$$ \end{enumerate}

\

\subsection{The Sol group}\label{sol}

\

The solvable group Sol topologically is the real vector space
${\mathbb{R} }^{3}$. Its Lie group structure is determined by the
product
\[(x,y,z)*(x',y',z')=(x+x',y+e^{-x}\,y',z+e^x\,z'), \]
where $(x,y,t),(x',y',t')\in \R^3$. See \cite{cannon}.

The left-invariant metric on Sol is

\bb\label{so0} ds^{2}=dx^2+e^{2x}\,dy^2+e^{-2x}\,dz^2 \ee and it
determines the semilinear Poisson equation

\bb\label{so1} u_{xx}+e^{-2x}u_{yy}+e^{2x}u_{zz}+f(u)=0. \ee

The sectional curvature of $(Sol,g)$ is nonconstant. See
\cite{bo}. Its scalar curvature $R=-2$.

The dimension of $Isom(Sol,g)$ is 3 (see \cite{bo,so}) and a basis
of Killing vector fields on Sol is given by \bb\label{so2}
So_{1}=\f{\p}{\p x}-y\f{\p}{\p y}+z\f{\p}{\p
z},\;\;\;\;\;So_{2}=\frac{\p}{\p y},\;\;\;\;\;So_{3}=\frac{\p}{\p
z}. \ee

\subsubsection{The group classification}

\begin{enumerate}
\item {\bf Arbitrary case:} The Lie point symmetry group of
(\ref{so1}) is $Isom(Sol,g)$ generated by (\ref{so2}). \item {\bf
Linear case:} In addition to the isometry group, we have the
symmetries
$$
So_{4}=u \frac{\p}{\p u},
$$
\bb\label{so3} So_{\infty}=b(x)\frac{\p}{\p u},\ee where $b$ is a
function such that $ \lb b +b=0.$ \item {\bf Homogeneous case}: We
have the same symmetries as in the linear case, but the function
$b$ in $(\ref{so3})$ satisfies $\lb b=0$ .
\end{enumerate}

\subsubsection{The Noether symmetries}

\begin{enumerate}
\item The isometry group $Isom(Sol,g)$ is a variational symmetry
group of (\ref{so1}) for any function $f(u)$. \item If $f(u)=0$ in
(\ref{so1}), the additional divergence symmetry is (\ref{so3}),
with $\lb b=0$. \item In the remaining linear case, the Noether
symmetries are $Isom(Sol)$ and (\ref{so3}), where $\lb b+b=0$.
\end{enumerate}

\subsubsection{The conservation laws}

\begin{enumerate}
\item For the symmetry $So_{1}$, with arbitrary $f(u)$, the
conservation law is $Div(A)=0$, where $A=(A_{1},A_{2},A_{3})$ and
$$
\ba{l c l}
A_{1} & = & \ds{\f{1}{2}(e^{-2x}u_{y}^{2}+e^{2x}u_{z}^{2}-u_{x}^{2})+yu_{x}u_{y}-zu_{x}u_{z}-F(u)},\\
\\
A_{2} & = & \ds{\f{1}{2}(yu_{x}^{2}+ye^{2x}u_{z}^{2}-ye^{-2x}u_{y}^{2})-e^{-2x}u_{x}u_{y}-e^{-2x}zu_{x}u_{z}+yF(u)},\\
\\
A_{3} & = & \ds{\f{z}{2}(u_{x}^{2}+e^{-2x}u_{y}^{2}+e^{2x}u_{z}^{2})-e^{2x}u_{x}u_{z}+e^{2x}yu_{y}u_{z}-e^{2x}zu_{y}u_{z}+yF(u)}.\\
\ea
$$

\item For the symmetry $So_{2}$, with arbitrary $f(u)$, the
conservation law is $Div(B)=0$, where $B=(B_{1},B_{2},B_{3})$ and
$$
\ba{l c l}
B_{1} & = & \ds{-u_{x}u_{y}},\\
\\
B_{2} & = & \ds{\f{1}{2}(u_{x}^{2}-e^{-2x}u_{y}^{2}+e^{2x}u_{z}^{2})-F(u)},\\
\\
B_{3} & = & \ds{-e^{2x}u_{y}u_{z}}.\\
\ea
$$

\item For the symmetry $So_{3}$, with arbitrary $f(u)$, the
conservation law is $Div(C)=0$, where $C=(C_{1},C_{2},C_{3})$ and
$$
\ba{l c l}
C_{1} & = & \ds{-u_{x}u_{z}},\\
\\
C_{2} & = & \ds{-e^{-2x}u_{y}u_{z}},\\
\\
C_{3} & = & \ds{\f{1}{2}(u_{x}^{2}+e^{-2x}u_{y}^{2}-e^{2x}u_{z}^{2})-F(u)}.\\
\ea
$$

\item For the symmetry $So_{\infty}$, with $\lb b=0$, the
conservation law is $Div(S)=0$, where $S=(S_{1},S_{2},S_{3})$ and
\bb\label{ss} \ba{l c l}
S_{1} & = & \ds{bu_{x}-b_{x}u},\\
\\
S_{2} & = & \ds{e^{-2x}bu_{y}-e^{-2x}b_{y}u},\\
\\
S_{3} & = & \ds{e^{2x}bu_{z}-e^{2x}b_{z}u}.\\
\ea \ee \item For the symmetry $So_{\infty}$, with $\lb b+b=0$,
the conservation law is $Div(S)=0$, where $S$ is given in
$(\ref{ss})$.
\end{enumerate}

\subsection{The product space $\mathbb{S}^{2}\times\R$}\label{prod1}

\

Let us now consider the set $\mathbb{S}^{2}\times\R$ endowed with
the product metric \bb\label{11.6.1}
ds^{2}=\f{dx^{2}+dy^{2}}{(1+x^{2}+y^{2})^{2}}+dz^{2}. \ee

The product manifold ($\mathbb{S}^{2}\times\R,g$) is a manifold
with constant scalar curvature $R=2$. The isometry group
$Isom(\mathbb{S}^{2}\times\R,g)$ has the following generators:
\bb\label{11.6.2} \ba{l c l}
S'_{1} & =& \ds{(1+x^{2}-y^{2})\f{\p}{\p x}+2xy\f{\p}{\p y}},\\
\\
S'_{2} &= & \ds{2xy\f{\p}{\p x}+(1-x^{2}+y^{2})\f{\p}{\p y}},\\
\\
S'_{3} &= & y\; \ds{\frac{\kd }{\kd x}-x\; \frac{\kd }{\kd y}},\\
\\
S'_{4} &= & \ds{\frac{\kd }{\kd z}} \ea \ee The nonlinear Poisson
equation on ($\mathbb{S}^{2}\times\R,g$) is given by
\bb\label{11.6.3}
(1+x^{2}+y^{2})^{2}(u_{xx}+u_{yy})+u_{zz}+f(u)=0. \ee
\subsubsection{Group classification}
\begin{enumerate}
\item {\bf Arbitrary case:} For any function $f(u)$, the symmetry
group coincides with $Isom(\mathbb{S}^{2}\times\R)$. \item {\bf
Homogeneous case:} If $f(u)=0$, then the additional symmetries are
\bb\label{11.6.4} S'_{\infty}=b\f{\p}{\p u}, \ee where $\lb b=0$
and \bb\label{11.6.5} S'_{5}=u\f{\p}{\p u}. \ee \item {\bf
Constant case:} The case $f(u)=k$ is reduced to the earlier under
the change $u\mapsto u-kz^{2}/2$. \item {\bf Linear case:} The
isometry group and the symmetry $S'_{\infty}$, with $\lb b=0$ in
(\ref{11.6.4}), generate a basis to the symmetry group generators.
\end{enumerate}

\subsubsection{The Noether symmetries}

\begin{enumerate}
\item The isometry group of $(\mathbb{S}^{2}\times\R,g)$ is a
variational symmetry group of equation (\ref{11.6.3}). \item
Symmetries (\ref{11.6.4}), with $\lb b=0$ or $\lb b+b=0$, are the
Noether symmetries in the cases $f(u)=0$ or $f(u)=u$,
respectively.
\end{enumerate}
\subsubsection{The conservation laws}
\begin{enumerate}
\item For the symmetry $S'_{1}$, with arbitrary $f(u)$, the
conservation law is $Div(A)=0$, where $A=(A_{1},A_{2},A_{3})$
and\newpage
$$
\ba{l c l}
A_{1} & = & \ds{\f{(1+x^{2}-y^{2})}{2}(u_{y}^{2}-u_{x}^{2})-2xyu_{x}u_{y}+\f{1+x^{2}-y^{2}}{2(1+x^{2}+y^{2})^{2}}u_{z}^{2}}\\
\\
&&\ds{-\f{1+x^{2}-y^{2}}{(1+x^{2}+y^{2})^{2}}F(u)},\\
\\
A_{2} & = & \ds{xyu_{x}^{2}-xyu_{y}^{2}-(1+x^{2}-y^{2})u_{x}u_{y}+\f{2xy}{(1+x^{2}+y^{2})^{2}}u_{z}^{2}}\\
\\
&&\ds{-\f{2xy}{(1+x^{2}+y^{2})^{2}}F(u)},\\
\\
A_{3} & = & \ds{-\f{1+x^{2}+y^{2}}{(1+x^{2}+y^{2})^{2}}u_{x}u_{z}-\f{xy}{(1+x^{2}+y^{2})^{2}}u_{y}u_{z}}.\\
\ea
$$

\item For the symmetry $S'_{2}$, with arbitrary $f(u)$, the
conservation law is $Div(B)=0$, where $B=(B_{1},B_{2},B_{3})$ and
$$
\ba{l c l}
B_{1} & = & \ds{xy(u_{y}^{2}-u_{x}^{2})+\f{xy}{(1+x^{2}+y^{2})^{2}}u_{z}^{2}-\f{2xy}{(1+x^{2}+y^{2})^{2}}F(u)},\\
\\
B_{2} & = & \ds{\f{1-x^{2}+y^{2}}{2}(u_{x}^{2}-u_{y}^{2})+\f{1-x^{2}+y^{2}}{2(1+x^{2}+y^{2})}u_{z}^{2}-2xyu_{x}u_{y}}\\
\\
&
&-\ds{\f{(1-x^{2}+y^{2})}{(1+x^{2}+y^{2})^{2}}F(u)},\\
\\
B_{3} & = & \ds{\f{2xy}{1+x^{2}+y^{2}}u_{x}u_{z}+\f{1-x^{2}+y^{2}}{(1+x^{2}+y^{2})^{2}}u_{y}u_{z}}.\\
\ea
$$

\item For the symmetry $S'_{3}$, with arbitrary $f(u)$, the
conservation law is $Div(C)=0$, where $C=(C_{1},C_{2},C_{3})$ and
$$
\ba{l c l}
C_{1} & = & \ds{\f{y}{2}(u_{y}^{2}-u_{x}^{2})+\f{y}{2(1+x^{2}+y^{2})^{2}}u_{z}^{2}+xu_{x}u_{z}-\f{y}{(1+x^{2}+y^{2})^{2}}F(u)},\\
\\
C_{2} & = & \ds{\f{x}{2}(u_{y}^{2}-u_{x}^{2})+\f{u_{z}^{2}}{(1+x^{2}+y^{2})^{2}}-yu_{x}u_{y}+\f{x}{(1+x^{2}+y^{2})^{2}}F(u)},\\
\\
C_{3} & = & \ds{\f{-y}{(1+x^{2}+y^{2})^{2}}u_{x}u_{y}+\f{-y}{(1+x^{2}+y^{2})^{2}}u_{y}u_{z}}.\\
\ea
$$

\item For the symmetry $S'_{4}$, with arbitrary $f(u)$, the
conservation law is $Div(D)=0$, where $D=(D_{1},D_{2},D_{3})$ and
$$
\ba{l c l}
D_{1} & = & -\ds{u_{x}u_{z}},\\
\\
D_{2} & = & -\ds{u_{y}u_{z},}\\
\\
D_{3} & = & \ds{\f{1}{2}(u_{x}^{2}+u_{y}^{2})-\f{u_{z}^{2}}{2(1+x^{2}+y^{2})^{2}}-\f{F(u)}{(1+x^{2}+y^{2})^{2}}}.\\
\ea
$$

\item For the symmetry $S'_{\infty}$, with $\lb b=0$ or $\lb
b+b=0$, the conservation law is $Div(E)=0$, where
$E=(E_{1},E_{2},E_{3})$ and \bb\label{sss1} \ba{l c l}
E_{1} & = & \ds{bu_{x}-b_{x}u},\\
\\
E_{2} & = & \ds{bu_{y}-b_{y}u},\\
\\
E_{3} & = & \ds{\f{bu_{z}-b_{z}u}{(1+x^{2}+y^{2})^{2}}}.\\
\ea \ee

\end{enumerate}

\

\subsection{The product space $\mathbb{H}^{2}\times\R$}\label{prod2}

\

Let $\mathbb{H}^{2} = \{(x, y) \in\R^{2}\, :\, y > 0\}$ endowed
with metric $$\f{dx^{2}+dy^{2}}{y^{2}}$$ be the hyperbolic plane
(Klein's Model) and consider the set $(H_{2}\times\R,g)$ endowed
with the product metric \bb\label{11.7.1}
ds^{2}=\f{dx^{2}+dy^{2}}{y^{2}}+dz^{2}. \ee

The Lie algebra of the infinitesimal isometries of
$(\mathbb{H}^{2}\times\R,g)$, $Isom(\mathbb{H}^{2}\times\R,g)$, is
given by (see \cite{irene}) \bb\label{11.7.2}
X_{1}=\f{x^{2}-y^{2}}{2}\f{\p}{\p x}+xy\f{\p}{\p
y},\,\,\,\,X_{2}=\f{\p}{\p x},\,\,\,\,X_{3}=x\f{\p}{\p
x}+y\f{\p}{\p y},\,\,\,\,X_{4}=\f{\p}{\p z}. \ee

The scalar curvature of $(\mathbb{H}^{2}\times\R,g)$ is $R=-2$ and
the nonlinear Poisson equation is given by \bb\label{11.7.3}
y^{2}(u_{xx}+u_{yy})+u_{zz}+f(u)=0. \ee

\subsubsection{Group classification}
\begin{enumerate}
\item {\bf Arbitrary case:} For any function $f(u)$, the symmetry
group coincides with $Isom(\mathbb{H}^{2}\times\R,g)$. \item {\bf
Homogeneous case:} If $f(u)=0$, then the additional symmetries are
\bb\label{11.7.4} X_{\infty}=b\f{\p}{\p u}, \ee where $\lb b=0$
and \bb\label{11.7.5} X_{5}=u\f{\p}{\p u}. \ee \item {\bf Constant
case:} The case $f(u)=k$ is reduced to the earlier under the
change $u\mapsto u-kz^{2}/2$. \item {\bf Linear case:} The
isometry group and the symmetry $X_{\infty}$, with $\lb b=0$ in
(\ref{11.7.4}), generate a basis to the symmetry group generators.
\end{enumerate}

\subsubsection{The Noether symmetries}
\begin{enumerate}
\item The isometry group of $Isom(\mathbb{H}^{2}\times\R,g)$ is a
variational symmetry group of equation (\ref{11.7.3}). \item
Symmetries (\ref{11.7.4}), with $\lb b=0$ or $\lb b+b=0$, are the
Noether symmetries in the cases $f(u)=0$ or $f(u)=u$,
respectively.
\end{enumerate}

\subsubsection{The conservation laws}

\begin{enumerate}
\item For the symmetry $X_{1}$, with arbitrary $f(u)$, the
conservation law is $Div(A)=0$, where $A=(A_{1},A_{2},A_{3})$ and
$$
\ba{l c l}
A_{1} & = & \ds{\f{x^{2}-y^{2}}{4}(u_{y}^{2}-u_{x}^{2})+\f{x^{2}-y^{2}}{4y^{2}}u_{z}^{2}-xyu_{x}u_{y}-\f{x^{2}-y^{2}}{2y^{2}}F(u)},\\
\\
A_{2} & = & \ds{\f{xy}{2}(u_{x}^{2}-u_{y}^{2})+\f{x}{2y}u_{z}^{2}-\f{x^{2}-y^{2}}{2}u_{x}u_{y}-\f{x}{y}F(u)},\\
\\
A_{3} & = & \ds{-\f{x^{2}-y^{2}}{2y^{2}}u_{x}u_{z}-\f{x}{y}u_{y}u_{z}}.\\
\ea
$$

\item For the symmetry $X_{2}$, with arbitrary $f(u)$, the
conservation law is $Div(B)=0$, where $B=(B_{1},B_{2},B_{3})$ and
$$
\ba{l c l}
B_{1} & = & \ds{\f{u_{y}^{2}-u_{x}^{2}}{2}+\f{u_{z}^{2}}{2y^{2}}-\f{F(u)}{y^{2}}},\\
\\
B_{2} & = & \ds{-u_{x}u_{y}}\\
\\
B_{3} & = & \ds{-u_{x}u_{z}}.\\
\ea
$$

\item For the symmetry $X_{3}$, with arbitrary $f(u)$, the
conservation law is $Div(C)=0$, where $C=(C_{1},C_{2},C_{3})$ and
$$
\ba{l c l}
C_{1} & = & \ds{\f{x}{2}(u_{y}^{2}-u_{x}^{2})+\f{x}{2y^{2}}u_{z}^{2}-yu_{x}u_{y}-\f{x}{y^{2}}F(u)},\\
\\
C_{2} & = & \ds{\f{y}{2}(u_{x}^{2}-u_{y}^{2})+\f{u_{z}^{2}}{2y^{2}}-xu_{x}u_{y}+\f{F(u)}{y}},\\
\\
C_{3} & = & \ds{-\f{x}{y^{2}}u_{x}u_{z}+\f{u_{y}u_{z}}{y}}.\\
\ea
$$

\item For the symmetry $X_{4}$, with arbitrary $f(u)$, the
conservation law is $Div(D)=0$, where $D=(D_{1},D_{2},D_{3})$ and
$$
\ba{l c l}
D_{1} & = & -\ds{u_{x}u_{z}},\\
\\
D_{2} & = & -\ds{u_{y}u_{z}}\\
\\
D_{3} & = & \ds{\f{u_{x}^{2}+u_{y}^{2}}{2}-\f{u_{z}^{2}}{2y^{2}}-\f{F(u)}{y^{2}}}.\\
\ea
$$

\item For the symmetry $X_{\infty}$, with $\lb b=0$ or $\lb
b+b=0$, the conservation law is $Div(E)=0$, where
$E=(E_{1},E_{2},E_{3})$ and \bb\label{sss2} \ba{l c l}
E_{1} & = & \ds{bu_{x}-b_{x}u},\\
\\
E_{2} & = & \ds{bu_{y}-b_{y}u},\\
\\
E_{3} & = & \ds{\f{bu_{z}-b_{z}u}{y^{2}}}.\\
\ea \ee
\end{enumerate}

\

\subsection{The universal covering of $SL_{2}(\R)$}\label{sl2}

\

The universal covering of the Lie group of $2\times 2$ matrices
with determinant equal to $1$, denoted by $\widetilde{SL_2}(\R)$,
topologically is $\R^{3}_{+}:=\{(x,y,z)\in\R^{3}\,:\,z>0\}$,
endowed with the Riemmanian metric \bb\label{11.8.1}
ds^{2}=\left(dx+\f{dy}{z}\right)^{2}+\f{dy^{2}+dz^{2}}{z^{2}}. \ee
$(\widetilde{SL_2}(\R),g)$ possesses scalar curvature $R=-5/2$.

The nonlinear Poisson equation induced by metric (\ref{11.8.1}) is
given by \bb\label{11.8.2}
2u_{xx}-2zu_{xy}+z^{2}(u_{yy}+u_{zz})+f(u)=0. \ee

We have not found references giving the explicit form of
$Isom(\widetilde{SL_2},g)$. (At least in the information sources
available to us.) However it is known that its dimension is 4. See
\cite{scott}. Then we shall proceed in a way opposite to that in
subsections 11.2-11.7.

Since the scalar curvature of $(\widetilde{SL_2},g)$ is constant
($R=-5/2$), from the proof of Theorem 1 we conclude that the
symmetry group of (\ref{11.8.2}), with $f(u)\neq \lambda
u,\,\lambda=const.$, is reduced to the symmetry group of the
arbitrary case. Then, using the package SYM \cite{grego1,grego2}
by Stelios Dimas et al., we obtain that the symmetries of equation
(\ref{11.8.2}) are determined by \bb\label{11.8.3} X_{1}=\f{\p}{\p
x},\,\,\,\,\,X_{2}=\f{\p}{\p y},\,\,\,\,\,\,X_{3}=y\f{\p}{\p
y}+z\f{\p}{\p z},\,\,\,\,\,\,X_{4}=z\f{\p}{\p
x}+\f{y^{2}-z^{2}}{2}+yz\f{\p}{\p z}. \ee

To see that (\ref{11.8.3}) is the isometry group of
$(\widetilde{SL_2},g)$ we have two alternatives. The first one is
to check whether the fields (\ref{11.8.3}) are Killing vector
fields. Indeed, a simple substitution of (\ref{11.8.3}) into the
Killing equations confirms this claim. Then, from \cite{scott},
$dim(Isom(\widetilde{SL_2},g))=4$. Thus, (\ref{11.8.3}) generate a
basis of the generators of the isometry group of
$(\widetilde{SL_2},g)$.

The second one is as follows. From Theorem 1 is proved that the
isometry group of $(\widetilde{SL_2},g)$ and the symmetry group of
(\ref{11.8.2}) are the same. Since (\ref{11.8.3}) generate a basis
of the symmetry group of (\ref{11.8.2}), we conclude that
(\ref{11.8.3}) is a basis of the isometry group of
$(\widetilde{SL_2},g)$.

This procedure suggests that the existing programs for symbolic
calculation of symmetries of differential equations may be used to
calculate the isometry group of the considered manifold.

\subsubsection{Group classification}
\begin{enumerate}
\item {\bf Arbitrary case:} For any function $f(u)$, the symmetry
group coincides with $Isom(\widetilde{SL_2}(\R),g)$. \item {\bf
Homogeneous case:} If $f(u)=0$, then the additional symmetries are
\bb\label{11.8.4} X_{\infty}=b\f{\p}{\p u}, \ee where $\lb b=0$
and \bb\label{11.8.5} X'=u\f{\p}{\p u}. \ee \item {\bf Constant
case:} The case $f(u)=k$ is reduced to the preceding case by the
change $u\mapsto u-kx^{2}/2$. \item {\bf Linear case:} The
isometry group and the symmetry $X_{\infty}$, with $\lb b+b=0$ in
(\ref{11.8.4}), generate a basis of the symmetry algebra.
\end{enumerate}

\subsubsection{The Noether symmetries}

\begin{enumerate}
\item The isometry group $Isom(\widetilde{SL_2}(\R),g)$ is a
variational symmetry group of equation (\ref{11.8.3}). \item
Symmetries (\ref{11.8.4}), with $\lb b=0$ or $\lb b+b=0$, are the
Noether symmetries to the cases $f(u)=0$ or $f(u)=u$,
respectively.
\end{enumerate}

\subsubsection{The conservation laws}

\begin{enumerate}
\item For the symmetry $X_{1}$, with arbitrary $f(u)$, the
conservation law is $Div(A)=0$, where $A=(A_{1},A_{2},A_{3})$ and
$$
\ba{l c l}
A_{1} & = & \ds{-\f{u_{x}^{2}}{z^{2}}+\f{u_{y}^{2}}{2}+\f{u_{y}^{2}}{2}-\f{F(u)}{z^{2}}},\\
\\
A_{2} & = & \ds{\f{u_{x}^{2}}{z}-u_{x}u_{y}},\\
\\
A_{3} & = & \ds{-u_{x}u_{z}}.\\
\ea
$$

\item For the symmetry $X_{2}$, with arbitrary $f(u)$, the
conservation law is $Div(B)=0$, where $B=(B_{1},B_{2},B_{3})$ and
$$
\ba{l c l}
B_{1} & = & \ds{-\f{2}{z^{2}}u_{x}u_{y}+\f{u_{y}^{2}}{z}},\\
\\
B_{2} & = & \ds{\f{u_{x}^{2}}{z^{2}}-\f{u_{y}^{2}}{2}+\f{u_{z}^{2}}{2}-\f{F(u)}{z^{2}}}\\
\\
B_{3} & = & \ds{-u_{y}u_{z}}.\\
\ea
$$

\item For the symmetry $X_{3}$, with arbitrary $f(u)$, the
conservation law is $Div(C)=0$, where $C=(C_{1},C_{2},C_{3})$ and
$$
\ba{l c l}
C_{1} & = & \ds{-\f{2y}{z}u_{x}u_{y}-\f{2}{z}u_{x}u_{z}+\f{y}{z}u_{y}^{2}+u_{y}u_{z}},\\
\\
C_{2} & = & \ds{\f{y}{z^{2}}u_{x}^{2}+u_{x}u_{y}-\f{y}{2}u_{y}^{2}-zu_{y}u_{z}+\f{y}{2}u_{z}^{2}-\f{y}{z^{2}}F(u)}\\
\\
C_{3} & = & \ds{\f{u_{x}^{2}}{z}-u_{x}u_{y}+\f{z}{2}u_{y}^{2}-yu_{y}u_{z}-\f{z}{2}u_{z}^{2}-\f{F(u)}{z}}.\\
\ea
$$

\item For the symmetry $X_{4}$, with arbitrary $f(u)$, the
conservation law is $Div(D)=0$, where $D=(D_{1},D_{2},D_{3})$ and
$$
\ba{l c l}
D_{1} & = & \ds{-\f{u_{x}^{2}}{z}+\f{z^{2}-y^{2}}{z^{2}}u_{x}u_{y}-\f{2y}{z}u_{x}u_{z}+\f{y^{2}}{2z}u_{y}^{2}+yu_{y}u_{z}+\f{z}{2}u_{z}^{2}-\f{F(u)}{z}},\\
\\
D_{2} & = & \ds{\f{y^{2}+z^{2}}{2z^{2}}u_{x}^{2}-zu_{x}u_{y}+yu_{x}u_{z}-yzu_{y}u_{z}+\f{y^{2}-z^{2}}{4}u_{z}^{2}+\f{z^{2}-y^{2}}{2z^{2}}F(u)}\\
\\
D_{3} & = & \ds{\f{y}{z}u_{x}^{2}-yu_{x}u_{y}-zu_{x}u_{z}+\f{yz}{2}u_{y}^{2}+\f{z^{2}-y^{2}}{2}u_{y}u_{z}-\f{yz}{2}u_{z}^{2}-\f{y}{z}F(u)}.\\
\ea
$$

\item For the symmetry $X_{\infty}$, with $\lb b=$ or $\lb b+b=0$,
the conservation law is $Div(E)=0$, where $E=(E_{1},E_{2},E_{3})$
and
$$
\ba{l c l}
E_{1} & = & \ds{\f{2}{z^{2}}(bu_{x}-b_{x}u)-\f{1}{z}(bu_{y}-b_{y}u)},\\
\\
E_{2} & = & \ds{-\f{1}{z}(bu_{x}-b_{x}u)+(bu_{y}-b_{y}u)}\\
\\
E_{3} & = & \ds{bu_{z}-b_{z}u}.\\
\ea
$$
\end{enumerate}

\

\subsection{The Heisenberg group}\label{heisenberg}

\

The three-dimensional nilpotent Lie group, also called Heisenberg
group $H^{1}$, topologically is the real three-dimensional vector
space $\R^{3}$ endowed with the group structure determined by the
composition law $\phi:\R^{3}\times\R^{3}\longrightarrow\R^{3}$
defined by
$$
\phi((x,y,t),(x',y',t')):=(x+x',y+y',t+t'+2(x'y-xy')).
$$

This composition law determines the left invariant vector fields
\bb\label{h1.1.9.1} T=\f{\p}{\p t},\,\,\,X=\f{\p}{\p
x}+2y\f{\p}{\p t},\,\,\,Y=\f{\p}{\p y}-2x\f{\p}{\p t} \ee and the
left invariant metric on $H^{1}$ \bb\label{h1.1.9.2}
ds^{2}=dx^2+dy^2+\big(dz+2ydx-2xdy\big)^2. \ee

The scalar curvature of $(H^{1},g)$ is $R=-8$ and the operators
(\ref{h1.1.9.1}) satisfy the following commutation relations:
$$[X,Y]=-4T,\,\,[X,T]=[Y,T]=0.$$

These formulae represent in an abstract form the commutation
relations for the quantum-mechanical position and momentum
operators. This justifies the name Heisenberg group.

It is well-known that the metric (\ref{h1.1.9.2}) determines the
following generators of the isometry group of $H^{1}$, denoted by
$Isom(H^{1},g)$: \bb\label{h1.9.1.3} T=\f{\p}{\p
t},\,\,\,\,\tilde{X}=\f{\p}{\p x}-2y\f{\p}{\p
t},\,\,\,\,\tilde{Y}=\f{\p}{\p y}+2x\f{\p}{\p
t},\,\,\,\,\,R=y\f{\p}{\p x}-x\f{\p}{\p y}. \ee Note that $T$
corresponds to translations in $t$, $R-$to rotations in the
$(x,y)$ plane and $ \tilde{X},\tilde{Y}$ determine the right
multiplication.

For more details, see \cite{yi1,yi2,irene}.

The nonlinear Poisson equation on $H^{1}$ is given by
\bb\label{h1.9.1.31}
u_{xx}+u_{yy}+[4(x^{2}+y^{2})+1]u_{tt}+4yu_{xt}-4xu_{yt}+f(u)=0.
\ee

\subsubsection{Group classification}
\begin{enumerate}
\item {\bf Arbitrary case:} For any function $f(u)$, the symmetry
group coincides with $Isom(H^{1},g)$. \item {\bf Homogeneous
case:} If $f(u)=0$, then the additional symmetries are
\bb\label{h1.1.9.3} H_{\infty}=b\f{\p}{\p u}, \ee where $\lb b=0$
and \bb\label{h1.1.9.4} H_{1}=u\f{\p}{\p u}. \ee \item {\bf
Constant case:} The case $f(u)=k$ is reduced to the earlier under
the change $u\mapsto u-kx^{2}$. \item {\bf Homogeneous case:} The
isometry group and the symmetry $H_{\infty}$, with $\lb b=0$ in
(\ref{h1.1.9.3}), generate a basis of the symmetry group
generators.
\end{enumerate}

\subsubsection{The Noether symmetries}

\begin{enumerate}
\item The isometry group of $(H^{1},g)$ is a variational symmetry
group of equation (\ref{h1.9.1.31}). \item Symmetries
(\ref{h1.1.9.3}), with $\lb b=0$ or $\lb b+b=0$, are the Noether
symmetries to the cases $f(u)=0$ or $f(u)=u$, respectively.
\end{enumerate}

\subsubsection{The conservation laws}

\begin{enumerate}
\item For the symmetry $T$, with arbitrary $f(u)$, the
conservation law is $Div(A)=0$, where $A=(A_{1},A_{2},A_{3})$ and
$$
\ba{l c l}
A_{1} & = & \ds{-u_{x}u_{t}-2yu_{t}^{2}},\\
\\
A_{2} & = & \ds{-u_{y}u_{t}+2xu_{t}^{2}},\\
\\
A_{3} & = & \ds{\f{u_{x}^{2}+u_{y}^{2}}{2}-\f{4(x^{2}+y^{2})+1}{2}u_{t}^{2}-F(u)}.\\
\ea
$$

\item For the symmetry $\tilde{X}$, with arbitrary $f(u)$, the
conservation law is $Div(B)=0$, where $B=(B_{1},B_{2},B_{3})$ and
$$
\ba{l c l}
B_{1} & = & \ds{\f{u_{y}^{2}-u_{x}^{2}}{2}++2yu_{x}u_{t}-2xu_{y}u_{t}+\f{4(x^{2}+3y^{2})+1}{2}u_{t}^{2}-F(u)},\\
\\
B_{2} & = & \ds{-u_{x}u_{y}-2yu_{y}u_{t}+2xu_{x}u_{t}-4xyu_{t}^{2}},\\
\\
B_{3} & = & \ds{-3yu_{x}^{2}-yu_{y}^{2}+2xu_{x}u_{y}+y[4(x^{2}+y^{2})+1]u_{t}^{2}}\\
\\
&&\ds{-[4(x^{2}+y^{2})+1]u_{x}u_{t}+2yF(u)}.\\
\ea
$$

\item For the symmetry $\tilde{Y}$, with arbitrary $f(u)$, the
conservation law is $Div(C)=0$, where $C=(C_{1},C_{2},C_{3})$ and
$$
\ba{l c l}
C_{1} & = & \ds{-u_{x}u_{y}-2xu_{x}u_{t}-2yu_{y}u_{t}-4xyu_{t}^{2}},\\
\\
C_{2} & = & \ds{\f{u_{x}^{2}-u_{y}^{2}}{2}+2yu_{x}u_{t}-2xu_{y}u_{t}+\f{4(3x^{2}+y^{2})+1}{2}u_{t}^{2}-F(u)},\\
\\
C_{3} & = & \ds{-xu_{x}^{2}+3xu_{y}^{2}-x[4(x^{2}+y^{2})+1]u_{t}^{2}-2yu_{x}u_{y}}\\
\\
&&\ds{-[4(x^{2}+y^{2})+1]u_{y}u_{t}-2xF(u)}.\\
\ea
$$

\item For the symmetry $R$, with arbitrary $f(u)$, the
conservation law is $Div(D)=0$, where $D=(D_{1},D_{2},D_{3})$ and
$$
\ba{l c l}
D_{1} & = & -\ds{\f{y}{2}(u_{y}^{2}-u_{x}^{2})+\f{y}{2}[4(x^{2}+y^{2})+1]u_{t}^{2}+xu_{x}u_{y}-yF(u)},\\
\\
D_{2} & = & -\ds{\f{x}{2}(u_{y}^{2}-u_{x}^{2})-\f{x}{2}[4(x^{2}+y^{2})+1]u_{t}^{2}-yu_{x}u_{y}+xF(u)}\\
\\
D_{3} & = & \ds{-2y^{2}u_{x}^{2}-2x^{2}u_{y}^{2}+4xyu_{x}u_{y}-y[4(x^{2}+y^{2})+1]u_{x}u_{t}}\\
\\
&&\ds{+x[4(x^{2}+y^{2})+1]u_{y}u_{t}}.\\
\ea
$$

\item For the symmetry $H_{\infty}$, with $\lb b=0$ or $\lb
b+b=0$, the conservation law is $Div(E)=0$, where
$E=(E_{1},E_{2},E_{3})$ and \bb\label{sss3} \ba{l c l}
E_{1} & = & \ds{b(u_{x}+2yu_{t})-u(b_{x}+2yu_{t})},\\
\\
E_{2} & = & \ds{b(u_{y}-2xu_{t})-u(b_{y}-2xu_{t})},\\
\\
E_{3} & = & \ds{b\{2yu_{x}-2xu_{y}+[4(x^{2}+y^{2})+1]u_{t}\}}\\
\\
&&\ds{-u\{2yb_{x}-2xb_{y}+[4(x^{2}+y^{2})+1]b_{t}\}}.\\
\ea \ee
\end{enumerate}

\

\section*{Acknowledgements}

 \

 \rm We wish to thank Professor Enzo Mitidieri for his firm encoragement.
 Yuri Bozhkov would also like to thank CNPq and FAPESP,
 Brasil, as well as FAEPEX-UNICAMP for partial
 financial support. Igor Leite Freire is thankful to IMECC-UNICAMP
 for gracious hospitality and UFABC for the leave which has
 given him the opportunity to visit UNICAMP where this paper was written.

\

\end{document}